\documentclass[12pt,twoside]{amsart}
\usepackage{amsfonts,amsmath}
\usepackage[all,knot,arc]{xy}

\def\endpf{\relax\ifmmode\expandafter\endproofmath\else
  \unskip\nobreak\hfil\penalty50\hskip.75em\hbox{}\nobreak\hfil\bull
  {\parfillskip=0pt \finalhyphendemerits=0 \bigbreak}\fi}
\def\bull{\vbox{\hrule\hbox{\vrule\kern3pt\vbox{\kern6pt}\kern3pt\vrule}\hrule}}

\newtheorem{defn}{Definition}[section]
\newtheorem{lemma}[defn]{Lemma}

\newtheorem{theorem}[defn]{Theorem}

\newtheorem{remark}[defn]{Remark}
\newtheorem{proposition}[defn]{Proposition}

\newtheorem{maintheorem}{Theorem}
\newtheorem{maincor}[maintheorem]{Corollary}

\newcommand{\cc}{{\mathbb C}}
\newcommand{\zz}{{\mathbb Z}}

\newcommand{\nn}{{\mathbb N}}
\newcommand{\qq}{{\mathbb Q}}
\newcommand{\pp}{{\mathbb P}}

\newcommand{\ozsvath}{Ozsv\'{a}th}
\newcommand{\szabo}{Szab\'{o}}
\newcommand{\saso}{Sa\v{s}o}
\newcommand{\spin}{\ifmmode{\rm Spin}\else{${\rm spin}$\ }\fi}
\newcommand{\spinc}{\ifmmode{{\rm Spin}^c}\else{${\rm spin}^c$\ }\fi}

\newcommand{\spincs}{\mathfrak s}

\newcommand{\lk}{{\rm lk}}

\newcommand{\calc}{\mathcal{C}}

\newcommand{\calq}{\mathcal{Q}}
\newcommand{\calu}{\mathcal{U}}

\newcommand{\ds}{\displaystyle}

\newenvironment{narrow}[2]{%
 \begin{list}{}{%
  \setlength{\topsep}{0pt}%
  \setlength{\leftmargin}{#1}%
  \setlength{\rightmargin}{#2}%
  \setlength{\listparindent}{\parindent}%
  \setlength{\itemindent}{\parindent}%
  \setlength{\parsep}{\parskip}%
 }%
\item[]}{\end{list}}

\newif\ifpic
\picfalse    
\pictrue   

\DeclareMathOperator\ch{Char}

\begin{document}

\title{Unknotting information from Heegaard Floer homology}
\author{Brendan Owens}
\date{\today}
\thanks{}

\begin{abstract}
We use Heegaard Floer homology to obtain bounds on unknotting
numbers. This is a generalisation of \ozsvath\ and \szabo's obstruction
to unknotting number one. We determine the unknotting numbers of
$9_{10}$, $9_{13}$, $9_{35}$, $9_{38}$, $10_{53}$, $10_{101}$ and
$10_{120}$; this completes the table of unknotting numbers for
prime knots with crossing number nine or less. Our obstruction
uses a refined version of Montesinos' theorem which gives a Dehn
surgery description of the branched double cover of a knot.
\end{abstract}

\maketitle

\pagestyle{myheadings} \markboth{BRENDAN OWENS}
{UNKNOTTING INFORMATION FROM HEEGAARD FLOER HOMOLOGY}


\section{Introduction}
\label{sec:intro}

Let $K$ be a knot in $S^3$.  Given any diagram $D$ for $K$, a new
knot may be obtained by changing one or more crossings of $D$. The
unknotting number $u(K)$ is the minimum number of crossing changes
required to obtain the unknot, where the minimum is taken over all
diagrams for $K$.

Let $\Sigma(K)$ denote the double cover of $S^3$ branched along
$K$. A theorem of Montesinos (\cite{mont}, or see Lemma
\ref{lem:monttrick}) tells us that for any knot $K$, $\Sigma(K)$
is given by Dehn surgery on some framed link in $S^3$ with $u(K)$
components, with half-integral framing coefficients. In particular
if $u(K)=1$ then $\Sigma(K)$ is obtained by $\pm \det K/2$ Dehn
surgery on a knot $C$, where $\det K$ is the determinant of $K$.
\ozsvath\ and \szabo\ have shown in \cite{osu1} that the Heegaard
Floer homology of a 3-manifold $Y$ gives an obstruction to $Y$
being given by half-integral surgery on a knot in $S^3$; they
apply this to $\Sigma(K)$ to obtain an obstruction to $K$ having
unknotting number one.

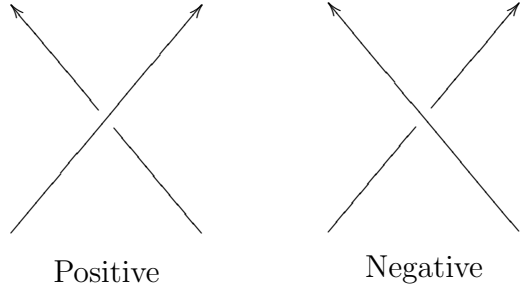
\begin{figure}[tbp] 
\begin{center}
\ifpic
\leavevmode
\begin{xy}
0;/r2pc/:
(0,0)*{
\begin{xy}
0;/r6pc/:
(0,0)*{}="1";
(1,0)*{}="2";
(1,1.2)*{}="3";
(0,1.2)*{}="4";
"1";"3" **\crv{}?(1)*\dir{>}; \POS?(.5)*{\hole}="x"; 
"2";"x" **\crv{}; 
"x";"4" **\crv{}?(1)*\dir{>}; 
(0.5,-0.2)*{{\rm Positive}}; 
\end{xy}};
(5,0)*{
\begin{xy}
0;/r6pc/:
(0,0)*{}="1";
(1,0)*{}="2";
(1,1.2)*{}="3";
(0,1.2)*{}="4";
"2";"4" **\crv{}?(1)*\dir{>}; \POS?(.5)*{\hole}="x"; 
"1";"x" **\crv{}; 
"x";"3" **\crv{}?(1)*\dir{>}; 
(0.5,-0.2)*{{\rm Negative}}; 
\end{xy}}
\end{xy}
\else \vskip 5cm \fi
\begin{narrow}{0.3in}{0.3in}
\caption{
\bf{Signed crossings in a knot diagram.}}
\label{fig:crossings}
\end{narrow}
\end{center}
\end{figure}

Note that crossings in a knot diagram may be given a sign as in
Figure \ref{fig:crossings} (independent of the choice of
orientation of the knot).  Let $\sigma(K)$ denote the signature of
a knot $K$. It is shown in \cite[Proposition 2.1]{cl} (also \cite[Theorem 5.1]{st}) 
that if $K'$ is obtained from
$K$ by changing a positive crossing, then
$$\sigma(K')\in\{\sigma(K),\sigma(K)+2\};$$
similarly if $K'$ is obtained from $K$ by changing a negative
crossing then
$$\sigma(K')\in\{\sigma(K),\sigma(K)-2\}.$$
Now suppose that $K$ may be unknotted by changing $p$ positive and
$n$ negative crossings (in some diagram).  Since the unknot has
zero signature, it follows that a bound for $n$ is given by
\begin{equation}
\label{eqn:nsig} n\ge\sigma(K)/2. \end{equation}

The main result of this paper is an obstruction to equality in
(\ref{eqn:nsig}).  This is easiest to state for the case of an
alternating knot; the obstruction is then a condition
on the Goeritz matrix obtained from an alternating
projection of $K$.  (We will
recall the definition of the Goeritz matrix in Section \ref{sec:os}.)
We also restrict for now to knots which can be unknotted with
two crossing changes.

A positive-definite integral matrix $Q$ of rank $r$
presents a finite group $\Gamma_Q$ via the short exact sequence
$$0\longrightarrow\zz^r\stackrel{Q}{\longrightarrow}\zz^r\longrightarrow\Gamma_Q\longrightarrow0.$$
A {\em characteristic covector} for $Q$ is an element of $\zz^r$
which is congruent modulo 2 to the diagonal of $Q$, i.e., an
element of
$$\ch(Q)=\{\left.\xi\in\zz^r\,\,\right|\,\,\xi_i\equiv Q_{ii}\pmod2\}.$$
Suppose that $\det Q$ is odd. Define a function
$$m_Q:\Gamma_Q\to\qq$$
by
$$m_Q(g) = \min\left\{\left.\frac{\xi^T Q^{-1}\xi-r}4\,\,\right|
\,\,\xi\in \ch(Q),\,[\xi]=g\right\}.$$
(The minimum exists since $Q$ is positive-definite.)

\begin{maintheorem}
\label{thm:u2}
Let $K$ be an alternating knot which may be unknotted by changing
$p$ positive and $n$ negative crossings, where $n=\sigma(K)/2$ and
$p+n=2$.  Let $G$ be the positive-definite Goeritz matrix obtained
from an alternating diagram for $K$.  Then there exists a positive-definite
matrix
$$
\widetilde{Q}=\left(\begin{matrix}
m_1 & 1 & a & 0\\
1 & 2 & 0 & 0\\
a & 0 & m_2 & 1\\
0 & 0 & 1 & 2
\end{matrix}\right),
$$
with
$$\det \widetilde{Q}=\det K,$$
$$0\le a<m_1\le m_2 \quad\mbox{(and hence $a<\det K/4$)},$$
and exactly $n$ of $\{m_1, m_2\}$ are even;
and a group isomorphism
$$\phi:\Gamma_{\widetilde{Q}}\to\Gamma_G$$
with
\begin{eqnarray*}
m_{\widetilde{Q}}(g)&\ge&m_G(g),\\
\mbox{and}\quad
m_{\widetilde{Q}}(g)&\equiv&m_G(g) \pmod2
\end{eqnarray*}
for all $g\in\Gamma_{\widetilde{Q}}$.
\end{maintheorem}

Applying Theorem \ref{thm:u2} to the alternating knots which were listed
in \cite{knotinfo} as having unknotting number 2 or 3 yields the
following:

\begin{maincor}
\label{cor:u} The knots $9_{10}, 9_{13}, 9_{35}, 9_{38}, 10_{53},
10_{101}, 10_{120}$ have unknotting number 3.
\end{maincor}

For all but one of the knots in Corollary \ref{cor:u}, the
signature is 4 and the unknotting number computation follows
directly from Theorem \ref{thm:u2}.
The exception is $9_{35}$, whose signature is 2.  The computation
of $u(9_{35})$ uses Theorem \ref{thm:u2} and also a result of
Traczyk \cite{t}.

Corollary \ref{cor:u} completes the table of unknotting numbers for prime
knots with 9 crossings or less.

Recall that for an oriented framed link $C_1,\ldots,C_r$ in $S^3$,
the linking matrix is the symmetric matrix $(a_{ij})$ with each
diagonal entry $a_{ii}$ given by the framing on $C_i$, and
off-diagonal entries $a_{ij}$ given by the linking numbers
$\lk(C_i,C_j)$. The following is a refinement of Montesinos'
theorem which was inspired by a theorem of Cochran and Lickorish \cite[Theorem 3.7]{cl}.

\begin{maintheorem}
\label{thm:mainthm} Suppose that a knot $K$ may be unknotted by
changing $p$ positive and $n$ negative crossings, with
$n=\sigma(K)/2$. Then the branched double cover
$\Sigma(K)$ may be obtained by Dehn surgery
on an oriented, framed $p+n$ component link
$C_1,\ldots,C_{p+n}$ in $S^3$
with linking matrix $\frac12Q$, where $Q$ is a positive-definite
integral matrix which is congruent to the identity modulo 2,
and exactly $n$ of the diagonal entries of $Q$ are congruent to
3 modulo 4.

Moreover, by handlesliding and changing orientations
one may replace the linking matrix with
$\frac12PQP^T$, for any $P\in GL(p+n,\zz)$ which is congruent to the identity
modulo 2.  This preserves the congruences modulo 4 on the diagonal.
\end{maintheorem}

It is shown in \cite{det3} that the double branched cover of the
Montesinos knot $10_{145}$ does not bound any positive-definite
four-manifold.  This knot has signature two.  Combining this with
Theorem \ref{thm:mainthm} (or the above-mentioned theorem of Cochran and Lickorish) 
yields the following:

\begin{maincor}
\label{cor:10_145}
If $10_{145}$ is unknotted by changing $p$ positive crossings and $n$
negative crossings, then $n\ge2$.
\end{maincor}

Given a matrix $Q$ in $M(r,\zz)$ which is conjugate modulo 2 to the identity, associate a matrix
$\widetilde{Q}\in M(2r,\zz)$ by replacing each entry by a $2\times2$-block as follows:
\begin{eqnarray}
\mbox{odd entries: }\qquad 2m-1 &\mapsto& \left[\begin{matrix}m&1\\1&2\end{matrix}\right]\nonumber\\
&\label{eqn:Q}\\
\mbox{even entries: }\qquad\quad 2a \quad &\mapsto&
\left[\begin{matrix}a&0\\0&0\end{matrix}\right].\nonumber
\end{eqnarray}

Thus for example if $r=2$,
$$Q=\left(\begin{matrix}2m_1-1&2a\\2a&2m_2-1\end{matrix}\right)\mapsto
\widetilde{Q}=\left(\begin{matrix}
m_1 & 1 & a & 0\\
1 & 2 & 0 & 0\\
a & 0 & m_2 & 1\\
0 & 0 & 1 & 2
\end{matrix}\right).
$$

For a rational homology three-sphere $Y$, the {\em correction terms} of
\ozsvath\ and \szabo\ are a set
of rational numbers $\{d(Y,\spincs)\,\,|\,\,\spincs\in\spinc(Y)\}$
which provide constraints on which four-manifolds $Y$ may bound.
We recall these constraints in Section \ref{sec:os}; combining these with
Theorem \ref{thm:mainthm}
yields the following unknotting obstruction, of which Theorem \ref{thm:u2} is a special
case.

\begin{maintheorem}
\label{thm:u}
Let $K$ be a knot in $S^3$ which may be unknotted by changing
$p$ positive and $n$ negative crossings, where $n=\sigma(K)/2$.
Let $Q_1,\ldots,Q_k$ be a complete set of representatives of
the finite quotient
$$\frac{\{Q\in M(p+n,\zz) \,\,|\,\,\mbox{\rm $Q$ is positive-definite, } \det Q = \det K,\,
Q\equiv I \pmod2\}}
{\{P\in GL(p+n,\zz)\,\,|\,\,P\equiv I \pmod2\}},$$
and let $\widetilde{Q}_1,\ldots,\widetilde{Q}_k$ be the corresponding elements of $M(2(p+n),\zz)$.
Then for some $Q_i$ which has exactly $n$ diagonal entries conjugate to 3 modulo 4, there exists
a group isomorphism
$$\phi:\Gamma_{\widetilde{Q}_i}\to\spinc(\Sigma(K))$$
with
\begin{eqnarray*}
m_{\widetilde{Q}_i}(g)&\ge&d(\Sigma(K),\phi(g)),\\
\mbox{and}\quad
m_{\widetilde{Q}_i}(g)&\equiv&d(\Sigma(K),\phi(g)) \pmod2
\end{eqnarray*}
for all $g\in\Gamma_{\widetilde{Q}_i}$.
\end{maintheorem}

The following example illustrates the use of Theorem \ref{thm:u}
to obstruct higher unknotting numbers.
\begin{maincor}
\label{cor:11a365}
The 11-crossing two-bridge knot $S(51,35)$
has unknotting number 4.
\end{maincor}

\vskip2mm \noindent{\bf Acknowledgements.}
The problem of generalising the obstruction in \cite{osu1} to higher unknotting numbers was
suggested to me by Peter \ozsvath.  I am grateful to Peter \ozsvath, Ravi Ramakrishna,
and \saso\ Strle for helpful discussions.
Some Maple programs used in verifying Corollaries \ref{cor:u} and \ref{cor:11a365}
were written jointly with \saso\ Strle.


\section{Kirby-Rolfsen calculus}
\label{sec:calc}

In this section we establish some preliminaries on Dehn surgery.  For details on
Dehn surgery and Kirby-Rolfsen calculus see \cite{gs}.

A framed link $L$ in $S^3$ with rational framing coefficients determines a three-manifold
$Y_L$ by Dehn surgery (remove a tubular neighbourhood of each component of $L$; the
framing coefficient determines the gluing map to sew back a solid torus along the boundary).
If the framing coefficients are integers one obtains a four-manifold $W_L$ with boundary
$Y_L$ by attaching two-handles to $B^4$ along the components of $L$.  Kirby-Rolfsen
calculus describes when two framed links $L,L'$ determine the same three-manifold $Y_L$.

Given a framed oriented link $L$ with components $C_1,\dots,C_m$,
let $A$ denote the free abelian
group with generators $c_1,\dots,c_m$.  Define a symmetric bilinear form
$$Q:A\times A\to\qq$$
by
$$Q(c_i,c_j)=\left\{
\begin{array}{ll}
\mbox{framing coefficient of $C_i$} & \mbox{if $i=j$;}\\
\mbox{linking number $\lk(C_i,C_j)$} & \mbox{if $i\ne j$.}
\end{array}\right.$$
In other words, the matrix of $Q$ in the basis $c_1,\dots,c_m$ is
the linking matrix of $L$.
(This is the intersection pairing on $H_2(W_L;\zz)$ if the diagonal entries are integers.)

In the case that the framing coefficients on $L$ are integers, any
change of basis in $A$ may be realised by a change in the link
$L$.  In particular the change of basis $c_i\mapsto c_i\pm c_j$
may be realised by a {\em handleslide}.  Let $\lambda_j$ denote a
pushoff of $C_j$ whose linking number with $C_j$ equals the
framing of $C_j$.  A handleslide $C_i\mapsto C_i\pm C_j$ consists
of replacing $C_i$ by the oriented band sum of $C_i$ with
$\pm\lambda_j$.  This gives a new link $L'$ whose linking matrix
is the matrix of $Q$ in the basis $c_1,\dots,c_i'=c_i\pm
c_j,\dots,c_m$ and with $Y_{L'}\cong Y_L$, $W_{L'}\cong W_L$. It
will be convenient to have the following generalisation of
handlesliding to links with rational framings.

\begin{proposition}
\label{prop:slide} Let $L$ be an oriented link in $S^3$ consisting
of components $C_1,\dots,C_m$ with framings
$\frac{p_1}{q_1},\dots,\frac{p_m}{q_m}$, and let $Q$ be the
rational-valued bilinear pairing determined by the linking matrix
of $L$. Then by replacing $C_i$ in $L$ it is possible to obtain a
link $L'$ whose linking matrix is the matrix of $Q$ in the basis
$c_1,\dots,c_i'=c_i\pm q_j c_j,\dots,c_m$ and with $Y_{L'}\cong
Y_L$.
\end{proposition}

\proof For each $j=1,\dots,m$ choose a continued fraction
expansion
$$\frac{p_j}{q_j} = a^j_{l_j}-\frac{1}{a^j_{l_j-1}-\raisebox{-3mm}{$\ddots$
\raisebox{-2mm}{${-\frac{1}{\displaystyle{a^j_1}}}$}}}\, .$$ (The
numbers $a^j_{l_j},\dots,a^j_1$ arise from the Euclidean algorithm
as follows:
\begin{eqnarray}
r_{l_j}=p_j&=&a^j_{l_j}q_j-r_{l_j-2}\nonumber\\
r_{l_j-1}=q_j&=&a^j_{l_j-1}r_{l_j-2}-r_{l_j-3}\nonumber\\
&\vdots&\label{eqn:Euc}\\
r_2&=&a^j_2r_1-1\nonumber\\
r_1&=&a^j_1\nonumber.)
\end{eqnarray}

\begin{figure}[htbp]
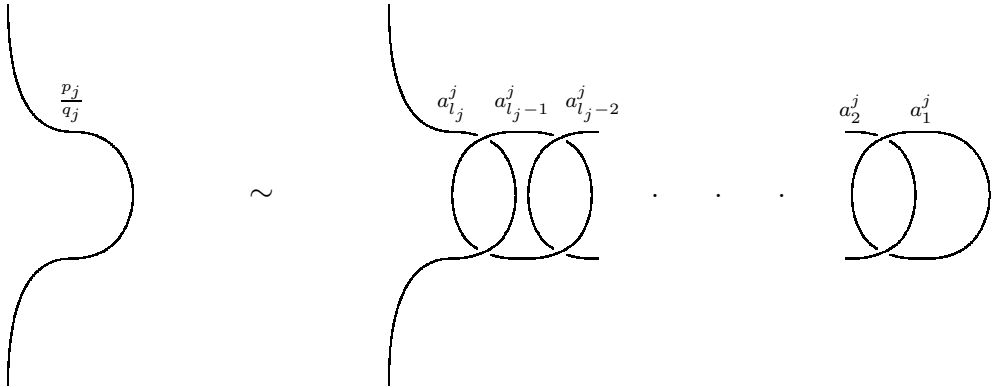
 
\begin{center}
\ifpic
\leavevmode
\xygraph{
!{0;/r2pc/:}
!{\xbendd[-2]@(0)} 
[d(3)l(1)]!{\xbendu[-2]@(0)} 
[u(2)r(0)]!{\xcaph[0.001]@(0)|{\frac{p_j}{q_j}}} 
[u(0)l(1)]!{\hcap[2]} 
[d(1)r(3)]!{*{\sim}}
[u(3)r(2)]!{\xbendd[-2]@(0)} 
[d(3)l(1)]!{\xbendu[-2]@(0)} 
[u(2)r(0)]!{\xcaph[0.001]@(0)|{a^j_{l_j}}} 
[u(0)l(1)]!{\vunder}!{\vunder-} 
[u(2)r(1)]!{\xcaph[0.2]@(0)|{a^j_{l_j-1}}} 
[d(2)l(1)]!{\xcaph[0.2]@(0)} 
[u(2)l(0.8)]!{\vunder}!{\vunder-} 
[u(2)r(1)]!{\xcaph[0.1]@(0)>{a^j_{l_j-2}}} 
[d(2)l(1)]!{\xcaph[0.1]@(0)} 
[u(1)r(0)]!{*{.}}
[r(1)]!{*{.}}
[r(1)]!{*{.}}
[u(1)r(1)]!{\xcaph[0.1]@(0)<{a^j_2}} 
[d(2)l(1)]!{\xcaph[0.1]@(0)} 
[u(2)l(0.9)]!{\vunder}!{\vunder-} 
[u(2)r(1)]!{\xcaph[0.2]@(0)|{a^j_1}} 
[d(2)l(1)]!{\xcaph[0.2]@(0)} 
[u(2)l(0.8)]!{\hcap[2]} 
}
\else \vskip 5cm \fi
\begin{narrow}{0.3in}{0.3in}
\caption{
\bf{Converting Dehn surgery to integral surgery.}}
\label{fig:qzsurg}
\end{narrow}
\end{center}
\end{figure}

Use reverse ``slam-dunks'' to obtain an integral surgery
description of $Y_L$: as shown in Figure \ref{fig:qzsurg}, we add
a chain of linked unknots linking each $C_j$, with framings
$a^j_1,\dots,a^j_{l_j-1}$, and replace the framing on $C_j$ with
$a^j_{l_j}$. (This is a standard procedure, see e.g. \cite[\S
5.3]{gs}.)  Denote the resulting link by $L_\zz$, and let
$Q_\zz:A_\zz\times A_\zz\to\zz$ denote the resulting bilinear
form.

\begin{figure}[tbp]
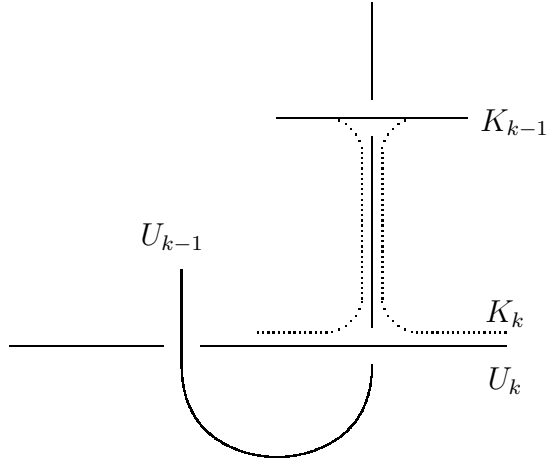
 
\begin{center}
\ifpic \leavevmode \xygraph{ !{0;/r6pc/:}
!{\xcaph[.8]@(0)}
[l(0.1)d(0.1)]!{\vcap[-1]}
[u(0.1)r(0.1)]!{\xcaph[1.6]@(0)}
[l(1.1)u(0.4)]!{\xcapv[0.5]@(0)}
[r(1)u(1.7)]!{\xcapv[1]@(0)}
[u(1.1)l(0.5)]!{\xcaph[1]@(0)}
[u(0.6)l(0.5)]!{\xcapv[0.5]@(0)}
[d(.3)l(1.575)]!{*\ellipse(0.2)_,=:a(80){.}}
[r(.25)d(0)]!{*\ellipse(0.2)_,=:a(80){.}}
!{(2.8475,1.015)*{}; (2.8475,0.25)*{}**\dir{.}}
!{(2.955,1.015)*{}; (2.955,0.25)*{}**\dir{.}}
[u(-.115)r(-1.38)]!{*\ellipse(0.2)_,=:a(80){.}}
[l(.25)u(0)]!{*\ellipse(0.2)_,=:a(80){.}}
!{(2.3,0.07)*{}; (2.66,0.07)*{}**\dir{.}}
!{(3.135,0.07)*{}; (3.6,0.07)*{}**\dir{.}}
[u(0.1)]!{*{K_k}} [d(0.35)]!{*{U_k}} [u(1.35)r(0.05)]!{*{K_{k-1}}}
[d(0.6)l(1.8)]!{*{U_{k-1}}} } \else \vskip 5cm \fi
\begin{narrow}{0.3in}{0.3in}
\caption{
\bf{Handlesliding $K_{k-1}$ over $U_k$ yields $K_k$ which is
separated from $U_{k-1}$ by a two-sphere.}} \label{fig:slide}
\end{narrow}
\end{center}
\end{figure}

We now perform handleslides on this integer-framed link.  Let
$U_1,\dots,U_{l_j-1}$ be the chain of unknots linking $C_j$ as
above, oriented so that $\lk(C_j,U_{l_j-1})=\lk(U_k,U_{k-1})=-1$,
for $2\le k<l_j$. Let $K_1=C_i+U_1$, and note that
\begin{eqnarray}
\label{eqn:K1} \lk(K_1,U_1)&=&a^j_1,\\
\nonumber \lk(K_1,U_2)&=&-1.
\end{eqnarray}
We now define $K_k$ recursively for $2\le k<l_j$. Choose any link
diagram of $K_{k-1}\cup U_{k-1}\cup U_k$. By performing a
handleslide over $U_k$ for each crossing where $K_{k-1}$ crosses
over $U_{k-1}$ we obtain a knot $K_k$ which does not cross over
$U_{k-1}$ and therefore is separated from it by a two-sphere in
$S^3$ (see Figure \ref{fig:slide}).  The signed count of these
handleslides is equal to the linking number of $K_{k-1}$ and
$U_{k-1}$; thus we write
$$[K_k] = [K_{k-1}]+\lk(K_{k-1},U_{k-1})[U_k],$$
where $[K_k]$ denotes the element of $A_\zz$ corresponding to the
knot $K_k$.  We may use this to
compute linking numbers and the framing of $K_k$.  In particular
\begin{equation}
\lk(K_2,U_2)=-1+a^j_2\,\lk(K_1,U_1),\label{eqn:K2}
\end{equation}
and for $2<k<l_j$,
\begin{eqnarray}
\lk(K_k,U_k)&=&\lk(K_{k-1},U_k)+a^j_k\,\lk(K_{k-1},U_{k-1})\label{eqn:Kk}\\
&=&-\lk(K_{k-2},U_{k-2})+a^j_k\,\lk(K_{k-1},U_{k-1})\nonumber .
\end{eqnarray}
Finally we let $C_i'$ be obtained as above from $K_{l_j-1}$ by
sliding over $C_j$, with $C_i'$ unlinked from each of
$U_1,\dots,U_{l_j-1}$.  We then have
\begin{eqnarray}
[C_i']&=&[K_{l_j-1}]+\lk(K_{l_j-1},U_{l_j-1})[C_j],\nonumber\\
\lk(C_i',C_j)&=&-\lk(K_{l_j-2},U_{l_j-2})+a^j_{l_j}\lk(K_{l_j-1},U_{l_j-1})
+\lk(C_i,C_j).\label{eqn:Ci'}
\end{eqnarray}

Comparing (\ref{eqn:K1}), (\ref{eqn:K2}), (\ref{eqn:Kk}), and
(\ref{eqn:Ci'}) to (\ref{eqn:Euc}) we see that
\begin{eqnarray*}
\lk(K_k,U_k)&=&r^j_k \quad\mbox{for $k=1,\dots,l_j-2$},\\
\lk(K_{l_j-1},U_{l_j-1})&=&r^j_{l_j-1}=q_j,\\
\lk(C_i',C_j)&=&p_j+\lk(C_i,C_j).
\end{eqnarray*}
This yields
$$[C_i']=[C_i]+\calu+q_j[C_j],$$
where
$$\calu=[U_1]+\ds\sum_{k=2}^{l_j-1}r_{k-1}[U_k].$$
Note that by construction $C_i'$ is separated by a two-sphere from
each $U_k$ and so $Q_\zz([C_i'],\calu)=0$. The framing of $C_i'$
is given by
\begin{eqnarray*}
Q_\zz([C_i'],[C_i'])&=&Q_\zz([C_i]+\calu+q_j[C_j],[C_i]+\calu+q_j[C_j])\\
&=&Q_\zz([C_i]+q_j[C_j],[C_i]+\calu+q_j[C_j])\\
&=&Q_\zz([C_i],[C_i])+2q_j Q_\zz([C_i],[C_j])+q_j^2a^j_{l_j}-q_jr_{l_j-2}\\
&=&a^i_{l_i}+2 q_j\lk(C_i,C_j) + p_j q_j.\\
\end{eqnarray*}
Slam dunking to remove the chains of linking unknots from each of
$C_1,\dots,C_i',\dots,C_m$ gives the required link $L'$ for the
basis change $c_i'=c_i+q_j c_j$. To get the opposite sign
construct $C_i'$ as above but start with $K_1=C_i-U_1$.\endpf

The following lemma is an application of the standard procedure,
referred to in the proof of Proposion \ref{prop:slide} and
illustrated in Figure \ref{fig:qzsurg}, for converting a Dehn
surgery description of a three-manifold to an integral surgery
description.

\begin{figure}[htbp]
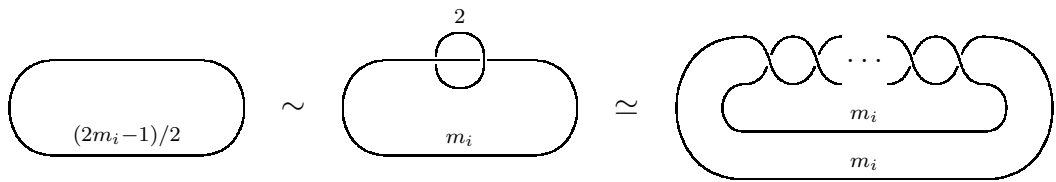
 
\begin{center}
\ifpic
\leavevmode
\xygraph{
!{0;/r1.5pc/:}
!{\hcap[-2]} 
!{\xcaph[3]@(0)}[d(2)l] 
!{\xcaph[3]@(0)|{(2m_i-1)/2}}[u(2)r(2)] 
!{\hcap[2]}[d(1)r(2)] 
!{*{\sim}}[u(1)r(2)]
!{\hcap[-2]} 
!{\xcaph[1.9]@(0)}[u(0.1)] 
!{\vcap|{2}}[d(0.2)] 
!{\vcap-}[u(0.2)r(1)] 
!{\xcapv[0.2]@(0)}[u(0.9)r(0.1)] 
!{\xcaph[0.9]@(0)}[d(2)l(3.1)] 
!{\xcaph[3]@(0)|{m_i}}[u(2)r(2)] 
!{\hcap[2]}[d(1)r(2)] 
!{*{\simeq}}[u(1.5)r(2.5)]
!{\hcap[-3]}[d] 
!{\hcap-}[u] 
!{\htwist}[ldd] 
!{\xcaph[1]@(0)}[ld] 
!{\xcaph[1]@(0)}[uuu] 
!{\htwist}[ldd] 
!{\xcaph[1]@(0)}[ld] 
!{\xcaph[1]@(0)}[u(2.5)r(0.5)] 
!{*{\dots}}[d(1.5)l(0.5)] 
!{\xcaph[1]@(0)|{m_i}}[ld] 
!{\xcaph[1]@(0)|{m_i}}[uuu] 
!{\htwist}[ldd] 
!{\xcaph[1]@(0)}[ld] 
!{\xcaph[1]@(0)}[uuu] 
!{\htwist}[ldd] 
!{\xcaph[1]@(0)}[ld] 
!{\xcaph[1]@(0)}[uuu] 
!{\hcap[3]}[d] 
!{\hcap} 
}
\else \vskip 5cm \fi
\begin{narrow}{0.3in}{0.3in}
\caption{ {\bf{Half-integer surgery.}} There are $2m_i-2$
crossings in the diagram on the right.} \label{fig:links}
\end{narrow}
\end{center}
\end{figure}

\begin{lemma}
\label{lem:links} Let $L=\{C_1,\ldots,C_r\}$ be a framed link in
$S^3$ with framing $(2m_i-1)/2$ on $C_i$, and let $Y$ be the
three-manifold obtained by Dehn surgery on $L$.  Then $Y$ is equal
to the boundary of the four-manifold $W$ obtained by adding
2-handles to $B^4$ along either of the following $2n$-component
framed links (as in Figure \ref{fig:links}):
\begin{enumerate}
\renewcommand{\labelenumi}{(\roman{enumi})}
\item the link consisting of the components $C_i$ with framing
$m_i$ plus a small linking unknot with framing $2$, for each
$i=1,\dots,r$;
\item the link consisting of $C_i$ with framing $m_i$, plus a
longitude $C'_i$ with framing $m_i$ and with the opposite orientation,
with linking number $\lk(C_i,C'_i)=1-m_i$, for each $i=1,\ldots,r$.
\end{enumerate}
\end{lemma}
\proof The fact that $Y$ is the boundary of the four-manifold
given by the framed link (i) follows from the continued fraction
expansions $(2m_i-1)/2=m_i-\frac12$.  The equivalence between (i)
and (ii) follows by handlesliding: add $C_i$ to $C'_i$ to go from
(ii) to (i).\endpf

Recall that to each matrix $Q\in M(r,\zz)$ which is congruent to the identity modulo 2, we
associate the matrix $\widetilde{Q}\in M(2r,\zz)$ as in (\ref{eqn:Q}).
If a 3-manifold $Y$ is given by Dehn surgery on a link with linking matrix $\frac12Q$,
then by Lemma \ref{lem:links}, $Y$ is the boundary of a simply-connected four-manifold
with intersection pairing $\widetilde{Q}$.  Also note that $\det Q=\det\widetilde{Q}$,
and $Q$ is positive-definite if and only if
$\widetilde{Q}$ is positive-definite: let $$\Delta_k(Q)=\det(Q_{ij})_{i,j\le k}.$$
Then
\begin{eqnarray*}
\Delta_{2k}(\widetilde{Q})&=&\Delta_k(Q),\\
\Delta_{2k-1}(\widetilde{Q})&=&(\Delta_{2k-2}(\widetilde{Q})+\Delta_{2k}(\widetilde{Q}))/2.
\end{eqnarray*}

\section{Proof of Theorem \ref{thm:mainthm}}
\label{sec:proof}

The proof of Theorem \ref{thm:mainthm} consists of three
lemmas.  The first of these is a proof of Montesinos' theorem
using Kirby calculus.  We could omit this and simply refer to
proofs in the literature, for example \cite{lick} (or to the
proof of Lemma \ref{lem:Q}).  We include the
proof since the four-dimensional point of view initially led us to
a proof of Theorem \ref{thm:mainthm}, and since it spells out a useful
algorithm for drawing a surgery diagram of $\Sigma(K)$.  
(For more details on Kirby diagrams of
cyclic branched covers see \cite[\S 6.3]{gs}; indeed what follows
is a variation of the method in their Exercise 6.3.5(c).)

\begin{lemma}
\label{lem:monttrick} Let $K$ be a knot in $S^3$ which can be
unknotted by changing $r$ crossings in some diagram $D$.  Then the
double branched cover $\Sigma(K)$ is given by Dehn surgery on an
$r$-component link in $S^3$ with linking matrix $\frac12Q$, where
$Q$ is congruent to the identity modulo 2.
\end{lemma}

\begin{figure}[htbp] 
\begin{center}
\ifpic
\leavevmode
\begin{xy}
(0,0)*{
\begin{xy}
0;/r4pc/:
(0.3,0.36)*{}="p1";
(1.7,2.04)*{}="p6";
(1.7,0.36)*{}="n1";
(0.3,2.04)*{}="n6";
"n1";"n6" **\crv{}?(1)*\dir{>}; \POS?(.5)*{\hole}="x"; 
"p1";"x" **\crv{}; 
"x";"p6" **\crv{}?(1)*\dir{>}; 
\end{xy}};
(30,-13)*{
\xygraph{
!{0;/r1.2pc/:}
[u(4)l(2)]!{\hcap[-2]} 
!{\xcaph[2]@(0)} 
[d(2)l]!{\xcaph[2]@(0)} 
[u(2)r(1)]!{\hcap[2]|{1}} 
[d(1)r(4)]!{*{\simeq}}
}};
(80,0)*{
\begin{xy}
0;/r4pc/:
(0.3,0.36)*{}="p1";
(.6,0.72)*{}="p2";
(.76,0.912)*{}="p3";
(.92,1.104)*{}="p4";
(1.08,1.296)*{}="p5";
(1.7,2.04)*{}="p6";
(1.7,0.36)*{}="n1";
(1.4,0.72)*{}="n2";
(1.24,0.912)*{}="n3";
(1.08,1.104)*{}="n4";
(.92,1.296)*{}="n5";
(0.3,2.04)*{}="n6";
"n1";"n2" **\crv{};
"n3";"n4" **\crv{};
"n5";"n6" **\crv{}?(1)*\dir{>};
"p1";"p2" **\crv{};
"p3";"p6" **\crv{}?(1)*\dir{>};
(1.4,1.68)*{}="parc";
(0.6,1.68)*{}="narc";
"parc";"narc" **\dir{.};
\end{xy}};
(80,-13)*{
\xygraph{
!{0;/r1.2pc/:}
[u(4)l(2)]!{\hcap[-2]} 
[d(2)]!{\xcaph[2]@(0)} 
[u(2)r(1)]!{\hcap[2]|{1}} 
[u(1.1)l(1)]!{*{\scriptstyle \alpha_i}} 
}};
(0,-40)*{
\begin{xy}
0;/r4pc/:
(0.3,0.36)*{}="p6";
(1.7,2.04)*{}="p1";
(1.7,0.36)*{}="n6";
(0.3,2.04)*{}="n1";
"n1";"n6" **\crv{}?(1)*\dir{>}; \POS?(.5)*{\hole}="x"; 
"p6";"x" **\crv{}; 
"x";"p1" **\crv{}?(1)*\dir{>}; 
\end{xy}};
(30,-53)*{
\xygraph{
!{0;/r1.2pc/:}
[u(4)l(2)]!{\hcap[-2]} 
!{\xcaph[2]@(0)} 
[d(2)l]!{\xcaph[2]@(0)} 
[u(2)r(1)]!{\hcap[2]|{1}} 
[d(1)r(4)]!{*{\simeq}}
}};
(80,-40)*{
\begin{xy}
0;/r4pc/:
(0.3,0.36)*{}="p6";
(.6,0.72)*{}="p5";
(.76,0.912)*{}="p4";
(.92,1.104)*{}="p3";
(1.08,1.296)*{}="p2";
(1.7,2.04)*{}="p1";
(1.7,0.36)*{}="n6";
(1.4,0.72)*{}="n5";
(1.24,0.912)*{}="n4";
(1.08,1.104)*{}="n3";
(.92,1.296)*{}="n2";
(0.3,2.04)*{}="n1";
"n1";"n2" **\crv{};
"n3";"n4" **\crv{};
"n5";"n6" **\crv{}?(1)*\dir{>};
"p6";"p5" **\crv{};
"p4";"p1" **\crv{}?(1)*\dir{>};
(1.4,1.68)*{}="parc";
(0.6,1.68)*{}="narc";
"parc";"narc" **\dir{.};
\end{xy}};
(80,-53)*{
\xygraph{
!{0;/r1.2pc/:}
[u(4)l(2)]!{\hcap[-2]} 
[d(2)]!{\xcaph[2]@(0)} 
[u(2)r(1)]!{\hcap[2]|{1}} 
[u(1.1)l(1)]!{*{\scriptstyle \alpha_i}} 
}}
\end{xy}
\else \vskip 5cm \fi
\begin{narrow}{0.3in}{0.3in}
\caption{
\bf{Changing crossings by sliding over a two-handle.}}
\label{fig:changes}
\end{narrow}
\end{center}
\end{figure}
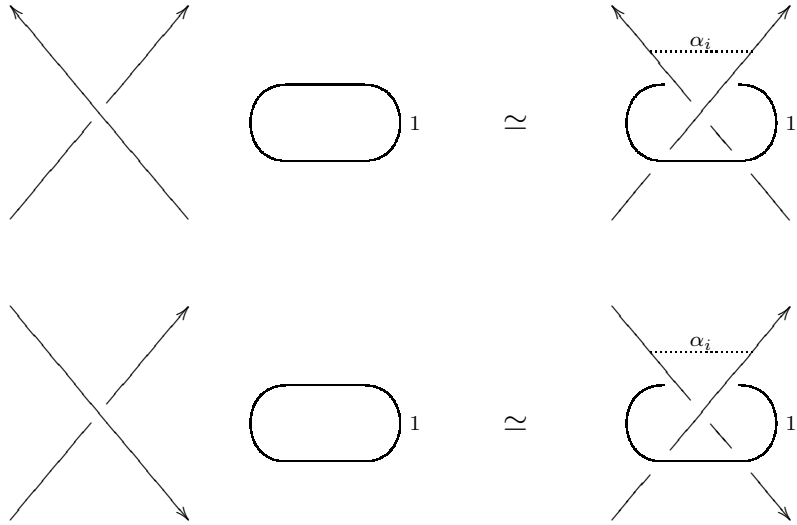

\proof We think of $K\subset S^3$ as being in the boundary of
$B^4$. Draw $r$ unlinked unknots beside $D$, each with framing
$+1$. This is a Kirby diagram which represents $K$ as a knot in
the boundary of the ``blown up'' four-ball $X=B^4\#\,r\cc\pp^2$.
As observed in \cite{cl}, the knot $K$ bounds a disk $\Delta$ in
$X$. This may be seen from the diagram by sliding each of the
chosen crossings in $D$ over a $+1$-framed unknot as in Figure
\ref{fig:changes}. Mark each of these changed crossings with a
small arc $\alpha_i$, $i=1,\dots,r$, as shown in that figure.

\begin{figure}[htbp]
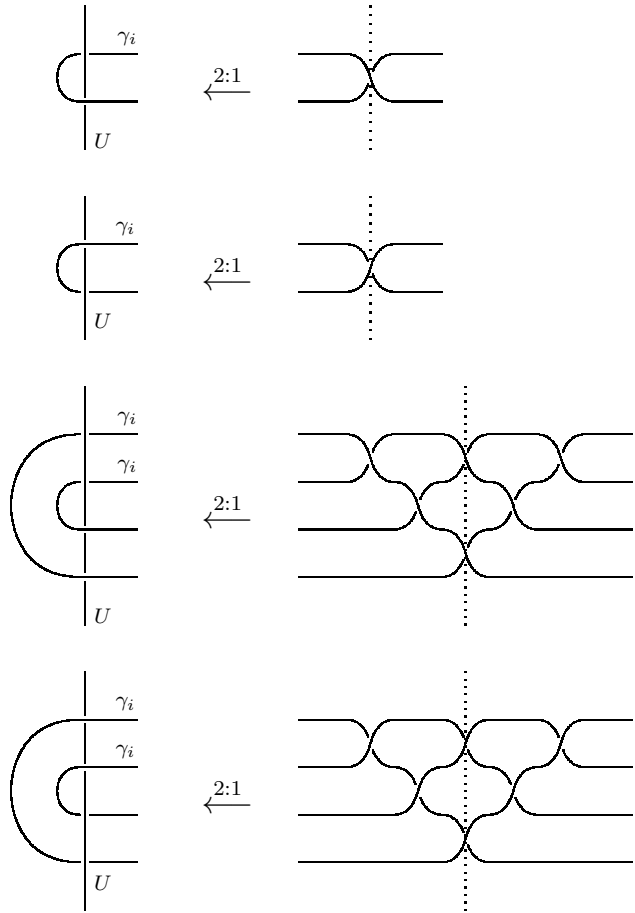
   
\begin{center}
\ifpic \leavevmode \xygraph{ !{0;/r1.5pc/:}
!{\hcap[-1]}[r(0.2)]
!{\xcaph[1]@(0)<{\gamma_i}}[d(1)l(1.2)] 
!{\xcaph[1.2]@(0)}[u(2)l(0.9)] 
!{\xcapv[1.9]@(0)}[d(1.1)] 
!{\xcapv[0.9]@(0)<{U}}[r(3)u(1.5)] 
!{*{\stackrel{2:1}{\longleftarrow}}}[r(1.5)u(.6)]
!{\xcaph[1]@(0)}[ld] 
!{\xcaph[1]@(0)}[u] 
!{\htwist} 
!{\xcaph[1]@(0)}[ld] 
!{\xcaph[1]@(0)}[u] 
!{\knotstyle{.}}[u(1)l(1.5)]
!{\xcapv[3]@(0)} 
!{\knotstyle{-}}[d(4)l(6.1)] !{\hcap[-1]}[r(0.2)d(1)]
!{\xcaph[1]@(0)}[u(1)l(1.2)] 
!{\xcaph[1.2]@(0)<{\gamma_i}}[u(1)l(0.9)] 
!{\xcapv[.9]@(0)}[d(0.1)] 
!{\xcapv[1.9]@(0)<{U}}[r(3)u(.5)] 
!{*{\stackrel{2:1}{\longleftarrow}}}[r(1.5)u(.6)]
!{\xcaph[1]@(0)}[ld] 
!{\xcaph[1]@(0)}[u] 
!{\htwistneg} 
!{\xcaph[1]@(0)}[ld] 
!{\xcaph[1]@(0)}[u] 
!{\knotstyle{.}}[u(1)l(1.5)]
!{\xcapv[3]@(0)}[d] 
!{\knotstyle{-}}[d(4)l(6.1)] !{\hcap[-1]}[r(0.2)]
!{\xcaph[1]@(0)<{\gamma_i}}[l(1.2)d(1)] 
!{\xcaph[1.2]@(0)}[u(2)l(1)] 
!{\hcap[-3]}[r(0.2)]
!{\xcaph[1]@(0)<{\gamma_i}}[l(1.2)d(3)] 
!{\xcaph[1.2]@(0)}[u(4)l(0.9)] 
!{\xcapv[2.9]@(0)}[d(2.1)] 
!{\xcapv[0.8]@(0)}[d(0)] 
!{\xcapv[0.9]@(0)<{U}}[r(3)u(2.5)] 
!{*{\stackrel{2:1}{\longleftarrow}}}[r(1.5)u(1.6)]
!{\xcaph[1]@(0)}[ld] 
!{\xcaph[1]@(0)}[ld] 
!{\xcaph[1]@(0)}[ld] 
!{\xcaph[1]@(0)}[u(3)] 
!{\htwist}[l(1)d(2)] 
!{\xcaph[1]@(0)}[ld] 
!{\xcaph[1]@(0)}[u(3)] 
!{\xcaph[1]@(0)}[ld] 
!{\htwist}[l(1)d(2)] 
!{\xcaph[1]@(0)}[u(3)] 
!{\htwist}[l(1)d(2)] 
!{\htwist}[u(2)] 
!{\xcaph[1]@(0)}[ld] 
!{\htwist}[l(1)d(2)] 
!{\xcaph[1]@(0)}[u(3)] 
!{\htwist}[l(1)d(2)] 
!{\xcaph[1]@(0)}[ld] 
!{\xcaph[1]@(0)}[u(3)] 
!{\xcaph[1]@(0)}[ld] 
!{\xcaph[1]@(0)}[ld] 
!{\xcaph[1]@(0)}[ld] 
!{\xcaph[1]@(0)} 
!{\knotstyle{.}}[u(4)l(3.5)]
!{\xcapv[5]@(0)}[d(3)l(2)] 
!{\knotstyle{-}}[d(4)l(6.1)] !{\hcap[-1]}[r(0.2)d(1)]
!{\xcaph[1]@(0)}[u(1)l(1.2)] 
!{\xcaph[1.2]@(0)<{\gamma_i}}[u(1)l(1)] 
!{\hcap[-3]}[r(0.2)d(3)]
!{\xcaph[1]@(0)}[u(3)l(1.2)] 
!{\xcaph[1.2]@(0)<{\gamma_i}}[u(1)l(0.9)] 
!{\xcapv[.9]@(0)}[d(0.1)] 
!{\xcapv[0.8]@(0)}[d(0)] 
!{\xcapv[2.9]@(0)<{U}}[r(3)u(0.5)] 
!{*{\stackrel{2:1}{\longleftarrow}}}[r(1.5)u(1.6)]
!{\xcaph[1]@(0)}[ld] 
!{\xcaph[1]@(0)}[ld] 
!{\xcaph[1]@(0)}[ld] 
!{\xcaph[1]@(0)}[u(3)] 
!{\htwistneg}[l(1)d(2)] 
!{\xcaph[1]@(0)}[ld] 
!{\xcaph[1]@(0)}[u(3)] 
!{\xcaph[1]@(0)}[ld] 
!{\htwistneg}[l(1)d(2)] 
!{\xcaph[1]@(0)}[u(3)] 
!{\htwistneg}[l(1)d(2)] 
!{\htwistneg}[u(2)] 
!{\xcaph[1]@(0)}[ld] 
!{\htwistneg}[l(1)d(2)] 
!{\xcaph[1]@(0)}[u(3)] 
!{\htwistneg}[l(1)d(2)] 
!{\xcaph[1]@(0)}[ld] 
!{\xcaph[1]@(0)}[u(3)] 
!{\xcaph[1]@(0)}[ld] 
!{\xcaph[1]@(0)}[ld] 
!{\xcaph[1]@(0)}[ld] 
!{\xcaph[1]@(0)} 
!{\knotstyle{.}}[u(4)l(3.5)]
!{\xcapv[5]@(0)} 
} \else \vskip 5cm \fi
\begin{narrow}{0.3in}{0.3in}
\caption{ {\bf{Drawing the double branched cover.}} Here $\Delta$
is the half-plane to the left of $U$.  The dotted line in the
diagrams on the right is the preimage of $U$.
The top two diagrams occur at endpoints of $\alpha_i$, and the bottom two occur
where $\alpha_i$ intersects the interior of $\Delta$.} \label{fig:2br1}
\end{narrow}
\end{center}
\end{figure}

\begin{figure}[htbp]
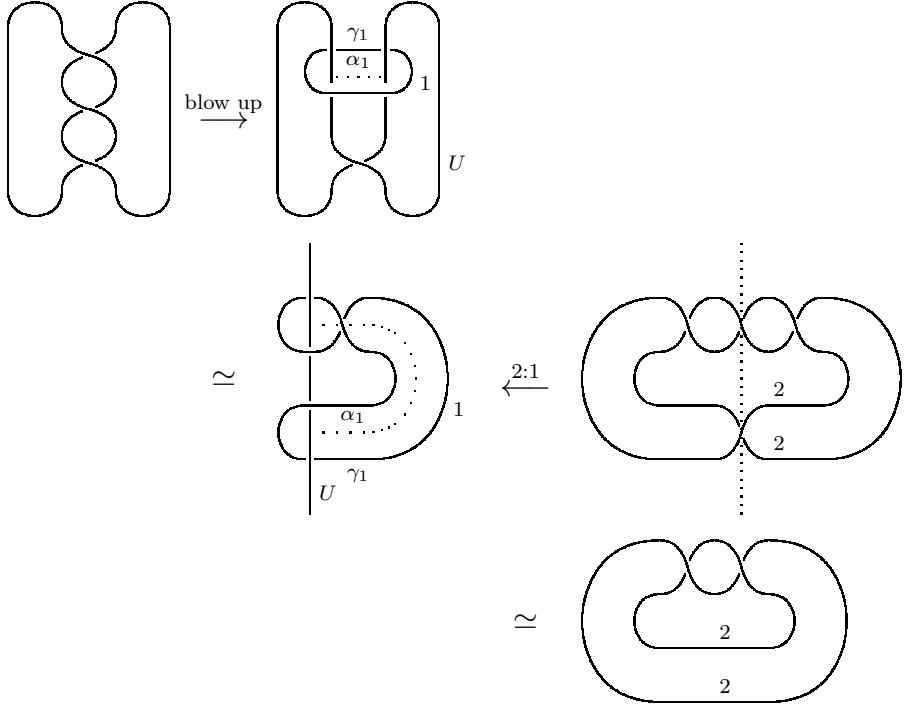
  
\begin{center}
\ifpic
\leavevmode
\xygraph{
!{0;/r1.7pc/:}
!{\vcap}[rr]
!{\vcap}[ll]
!{\xcapv[1]@(0)}[ur]
!{\vtwistneg}[urr]
!{\xcapv[1]@(0)}[lll]
!{\xcapv[1]@(0)}[ur]
!{\vtwistneg}[urr]
!{\xcapv[1]@(0)}[lll]
!{\xcapv[1]@(0)}[ur]
!{\vtwistneg}[urr]
!{\xcapv[1]@(0)}[lll]
!{\vcap-}[rr]
!{\vcap-}[u(1.5)r(2)]
!{*{\stackrel{\mbox{\tiny blow up}}{\longrightarrow}}}[u(1.5)r]
!{\vcap}[rr]
!{\vcap}[ll]
!{\xcapv[2]@(0)}[ur]
!{\xcapv[1]@(0)}[ur]
!{\xcapv[1]@(0)}[ur]
!{\xcapv[2]@(0)}[l(2.1)u(0.6)]
!{\hcap[-0.8]}[d(0.8)]
!{\xcaph[1.2]@(0)}[u(0.8)r(0.2)]
!{\hcap[0.8]<{1}}[l]
!{\xcaph[0.8]@(0)|{\gamma_1}}[d(0.5)l]
!{\knotstyle{.}\xcaph[0.8]@(0)|{\alpha_1}}[l(1.1)d(0.4)]
!{\knotstyle{-}\xcapv[0.7]@(0)}[ur]
!{\xcapv[0.7]@(0)}[u(0.3)l(2)]
!{\xcapv[1]@(0)}[ur]
!{\vtwistneg}[urr]
!{\xcapv[1]@(0)|{U}}[lll]
!{\vcap-}[rr]
!{\vcap-}[d(3.5)l(3)]
!{*{\simeq}}[u(1.5)r(1.5)]
!{\hcap[-1]}[d(2)]
!{\hcap[-1]}[u]
!{\xcaph[0.2]@(0)}[dl]
!{\xcaph[0.2]@(0)}[u(2)l(0.8)]
!{\htwist}[ldd]
!{\xcaph[1]@(0)}[ld]
!{\xcaph[-1]@(0)>{\gamma_1}}[uuu]
!{\hcap[3]<{1}}[d]
!{\hcap}[u(2)l(1.1)]
!{\xcapv[1.9]@(0)}[d(1.1)]
!{\xcapv[0.8]@(0)}[d(0)]
!{\xcapv[1.9]@(0)<{U}}[u(2.6)]
!{\knotstyle{.}}
!{\xcaph[1]@(0)}[d(2)l]
!{\xcaph[1]@(0)<{\alpha_1}}[u(2)]
!{\hcap[2]}[d(1)r(3)]
!{*{\stackrel{2:1}{\longleftarrow}}}[r(2.5)u(1.5)]
!{\knotstyle{-}}
!{\hcap[-3]}[d] 
!{\hcap-}[u] 
!{\htwist}[ldd] 
!{\xcaph[1]@(0)}[ld] 
!{\xcaph[1]@(0)}[uuu] 
!{\htwist}[ldd] 
!{\htwistneg}[uu] 
!{\htwist}[ldd] 
!{\xcaph[1]@(0)>{2}}[ld] 
!{\xcaph[1]@(0)>{2}}[uuu] 
!{\hcap[3]}[d] 
!{\hcap} 
!{\knotstyle{.}}[u(2)l(1.5)]
!{\xcapv[5]@(0)}[d(6)l(4)] 
!{*{\simeq}}[r(2.5)u(1.5)]
!{\knotstyle{-}}
!{\hcap[-3]}[d] 
!{\hcap-}[u] 
!{\htwist}[ldd] 
!{\xcaph[1]@(0)}[ld] 
!{\xcaph[1]@(0)}[uuu] 
!{\htwist}[ldd] 
!{\xcaph[1]@(0)>{2}}[ld] 
!{\xcaph[1]@(0)>{2}}[uuu] 
!{\hcap[3]}[d] 
!{\hcap} 
}
\else \vskip 5cm \fi
\begin{narrow}{0.3in}{0.3in}
\caption{
{\bf{A four-manifold bounded by the double branched cover of the left-handed trefoil.}}
}
\label{fig:2br2}
\end{narrow}
\end{center}
\end{figure}

The resulting diagram consists of:
\begin{itemize}
\item an unknot $U$ which has been obtained from $K$ by crossing changes;
\item arcs $\alpha_1,\dots,\alpha_r$ (one per changed crossing);
\item $+1$-framed unknots $\gamma_1,\dots,\gamma_r$.
\end{itemize}
Each $\gamma_i$ bounds a disk $D_i$ which retracts onto $\alpha_i$
and whose intersection with $U$ consists of the endpoints of
$\alpha_i$.

It is also observed in \cite{cl} that $H_1(X-\Delta;\zz/2)\cong\zz/2$,
with generator given by the meridian of $K$.  (To see this note
from Figure \ref{fig:changes} that the linking number of $U$ with each
of the $+1$-framed unknots is even.  Now use the Mayer-Vietoris sequence
for the decomposition of $X$ into $X-\Delta$ and a neighbourhood of
$\Delta$, with $\zz/2$ coefficients.)
Thus there exists a unique double cover $W$
of $X$ branched along $\Delta$; this is a four-manifold
with boundary $\Sigma(K)$.

Rearrange the diagram so that a point of $U$ which is not the
endpoint of an arc $\alpha_i$ is the point at infinity and $U$ is
a vertical line; then $\Delta$ may be seen in this diagram as the
half-plane to the left of $U$. (For a simple example see the first 3 diagrams in
Figure \ref{fig:2br2}.  Note in general the arcs $\alpha_i$ may be knotted
and linked, and may intersect $\Delta$.)  We may rearrange the diagram so that all intersections of $\gamma_i$
and $\Delta$ look like one of the diagrams on the left of Figure \ref{fig:2br1}.
To draw a Kirby diagram of $W$, we simply need to take two copies of $S^3-U$ cut
open along $\Delta$, and join the boundary half-planes in pairs.
Or in other words: take the part of the diagram
to the right of $U$, and draw another copy of it to the left of $U$. (Think of rotating the
half plane to the right of $U$ about $U$ by $\pi$, \emph{not} reflecting.)
Complete the centre of the diagram using Figure \ref{fig:2br1}.  (For an example see
Figure \ref{fig:2br2}.)

Each arc $\alpha_i$ lifts to a knot $\tilde{\alpha_i}$, and each
$D_i$ lifts to an annulus $\tilde{D_i}$ with core
$\tilde{\alpha_i}$.  The knot $\gamma_i$ lifts to two knots $C_i$,
$C'_i$; these are the boundary of the annulus $\tilde{D_i}$.

We now compute the framings of $C_i,C'_i$.  The $0$-framing on $\gamma_i$
lies on the disk $D_i$, and lifts to a curve on the annulus $\tilde{D_i}$.
This is the same framing for $C_i$ (or $C'_i$) as the other boundary curve
of $\tilde{D_i}$, but with the opposite sign.  Thus the $0$-framing on $\gamma_i$
lifts to the $-\lk(C_i,C'_i)$-framing on each of $C_i$, $C'_i$.  
Then the framing $+1$ on $\gamma_i$ lifts to $m_i$ on each of
$C_i,C'_i$, where $\lk(C_i,C'_i)=1-m_i$.

We note that the resulting Kirby diagram for $W$ matches that in
Lemma \ref{lem:links} (ii).  That lemma then shows that
$\Sigma(K)=\partial W$ is Dehn surgery on the framed link $L =
C_1,\dots,C_r$ with framing $(2m_i-1)/2$ on $C_i$.\endpf

To prove that the matrix $Q$ is positive-definite under the
hypotheses of Theorem \ref{thm:mainthm} one may appeal to
\cite[Theorem 3.7]{cl}, which gives a formula for the signature of
the four-manifold $W$ constructed in Lemma \ref{lem:monttrick}.
Surprisingly however it is also possible to prove this using the
following purely three-dimensional argument.

\begin{lemma}
\label{lem:Q} Suppose that $K$ may be unknotted by changing $p$
positive and $n$ negative crossings, with $n=\sigma(K)/2$.  Let
$\Sigma(K)$ be given by Dehn surgery on a link
$C_1,\dots,C_{p+n}$ with linking matrix $\frac12Q$ as in Lemma
\ref{lem:monttrick}.  Then $Q$ is positive-definite, and exactly
$n$ of the diagonal entries of $Q$ are congruent to 3 modulo 4.
\end{lemma}

\proof  The positivity of $Q$ is proved in \cite[Theorem 8.1]{osu1} for the case of
unknotting number one knots, i.e. $p+n=1$.  We include the proof of this case here
for completeness.  

\begin{figure}[htbp]   
\begin{center}
\ifpic \leavevmode \xygraph{ !{0;/r4.0pc/:}
!{\xunderv-=>}
[dr(0.5)]!{*{K_-}}
[ur(1.5)]!{\xunoverv-=>}
[dr(0.5)]!{*{K_0}}
[ur(1.5)]!{\xoverv-=>}
[dr(0.5)]!{*{K_+}}
} \else \vskip 5cm \fi
\begin{narrow}{0.3in}{0.3in}
\caption{} 
\label{fig:skein}
\end{narrow}
\end{center}
\end{figure}

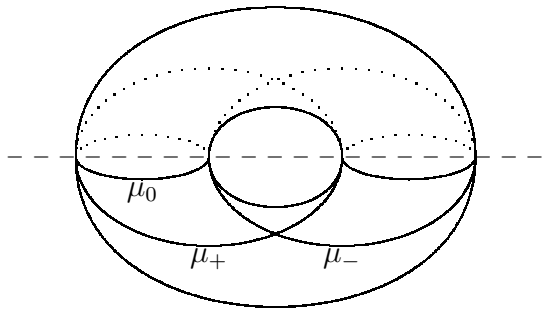
\begin{figure}[htbp] 
\begin{center}
\ifpic
\leavevmode
\begin{xy}
0;/r0.7pc/:
(-3,0)*\ellipse(3,1){.};
(-3,0)*\ellipse(3,1)__,=:a(-180){-};
(3,0)*\ellipse(3,1){.};
(3,0)*\ellipse(3,1)__,=:a(-180){-};
(-1.5,0)*\ellipse(6,4){.};
(-1.5,0)*\ellipse(6,4)__,=:a(-180){-};
(1.5,0)*\ellipse(6,4){.};
(1.5,0)*\ellipse(6,4)__,=:a(-180){-};
(-9,0)*{}="X1";
(-3,0)*{}="X2";
(3,0)*{}="X3";
(9,0)*{}="X4";
"X1";"X4" **\crv{(-9,9) & (9,9)};
"X1";"X4" **\crv{(-9,-9) & (9,-9)};
"X2";"X3" **\crv{(-3,3) & (3,3)};
"X2";"X3" **\crv{(-3,-3) & (3,-3)};
(-12,0)*{};(12,0)*{}**\dir{--};
(-6,-1.6)*{\mu_0};
(-3,-4.6)*{\mu_+};
(3,-4.6)*{\mu_-};
\end{xy}
\else \vskip 5cm \fi
\begin{narrow}{0.3in}{0.3in}
\caption{
{\bf{Meridians in $\Sigma(K_-)$, $\Sigma(K_0)$, $\Sigma(K_+)$.}}
Rotating by $\pi$ around the horizontal axis gives the solid torus as the
double cover of a ball branched along the arcs in $K_0$.}
\label{fig:torus}
\end{narrow}
\end{center}
\end{figure}

Suppose $K_-$, $K_0$ and $K_+$ are links in $S^3$ which are
identical outside of a ball in which they appear as in Figure \ref{fig:skein}.
Recall that the double cover of a ball $B$ branched along two arcs is a 
solid torus $\tilde B$, and
a meridian for the solid torus is given by the preimage in $\tilde B$ of either of the arcs
pushed out to the boundary of $B$.  It follows that $\Sigma(K_-)$, $\Sigma(K_0)$, $\Sigma(K_+)$
each contain an embedded solid torus, such that the complements of these solid tori can be identified.
The meridians which bound in  $\Sigma(K_-)$, $\Sigma(K_0)$, $\Sigma(K_+)$ are shown in Figure
\ref{fig:torus}.  They may be oriented so that their homology classes intersect as follows:
\begin{equation}
\label{eqn:int}
\mu_-\cdot \mu_+=2,\ \mu_+\cdot \mu_0=\mu_-\cdot \mu_0=1.
\end{equation}

Suppose now that $K=K_-$ with $\sigma(K)=2$, and $K_+$ is the unknot.  Then $\Sigma(K_+)=S^3$,
and $\Sigma(K_-)$ is $(2m-1)/2$ surgery on some knot $C$.  
We wish to show that $m$ is positive and even.
For some longitude $\lambda$
of $C$ with $\mu_+\cdot\lambda=1$ we have
$$\mu_-=-2\lambda-(2m-1)\mu_+.$$
Expressing $\lambda$ in the basis $\mu_+$, $\mu_0$ and plugging into (\ref{eqn:int}) yields 
$\lambda=\mu_0-m\mu_+$, from which we see that $\mu_0=\lambda+m\mu_+$.
In other words, $\Sigma(K_0)$ is $m$ surgery on $C$.

We now use two properties of the Conway-normalised Alexander polynomial, c.f. \cite{lick2}.
Firstly, for a knot $K$, the sign of the Alexander polynomial at $-1$ is given by 
$$(-1)^{\sigma(K)/2}\det K=\Delta_K(-1).$$  
This shows that $\Delta_{K_+}(-1)=1$ and $\Delta_{K_-}(-1)=-|2m-1|$.
Secondly we have the skein relation
$$\Delta_{K_+}(t)-\Delta_{K_-}(t)=(t^{-1/2}-t^{1/2})\Delta_{K_0}(t),$$
which yields
$$1+|2m-1|=2|m|.$$
It follows that $m$ and $2m-1$ are both positive.  Finally, the determinant
and signature of a knot $K$ are
shown in \cite[Theorem 5.6]{m} to satisfy
\begin{equation}
\label{eqn:detsig}
\det(K)\equiv\sigma(K)+1\pmod4,
\end{equation}
from which it follows that $2m-1$ is congruent to 3 modulo 4.

Similarly if $K_-$ is the unknot and $\sigma(K_+)=0$, we have that $\Sigma(K_+)$
is $(2m-1)/2$ surgery on a knot $C$ and we find $\Sigma(K_0)$ is $(m-1)$ surgery on $C$.
The skein relation gives $|2m-1|-1=2|m-1|$, which again shows $m$ is positive.
From (\ref{eqn:detsig}) we have $2m-1$ is congruent to 1 modulo 4.

The general case follows easily from the above.
Let $c_1,\dots,c_{p+n}$ be the set of crossings ($p$ positive, $n$ negative) 
in some chosen diagram that we change to unknot $K$.
Then $\Sigma(K)$ is Dehn surgery on a link
$L=C_1,\dots,C_{p+n}$, with linking matrix $\frac12Q$.  Each $C_i$
corresponds to a crossing $c_i$.
Dehn surgery on a sublink of $L$ gives the double
branched cover of a knot which is obtained from $K$ by changing
a subset of the crossings in $\calc$.  In particular $Q_{ii}/2$ surgery on
the knot $C_i$ yields the double branched cover of the knot
$K'$ which is obtained from $K$ by changing all of the crossings except $c_i$.  
It follows from the unknotting number one case applied to
$K'$ that all diagonal entries of $Q$ are positive and exactly those which 
correspond to negative crossings are congruent to 3 modulo 4.

It remains to prove that $Q$ is positive-definite.  
Note that from (\ref{eqn:nsig}) and the assumption $n=\sigma(K)/2$, the knot signature changes every
time we change a negative crossing and is unchanged when we change a positive crossing.
Let $Q_k$ be the submatrix
$\left(Q_{ij}\right)_{i,j\le k}$.  Observe that since the off-diagonal entries are even, 
the determinant of $Q_k$ is congruent modulo 4 to the product of the diagonal entries.  
Let $K_k$ be the knot obtained from $K$ by changing the crossings 
$c_{k+1},\dots,c_{p+n}$.  Suppose that $\det Q_{k-1}$ is positive, and hence equals
$\det K_{k-1}$.  If $c_k$ is positive then 
$$Q_{kk}\equiv 1 \implies \det Q_k \equiv \det Q_{k-1}\pmod4.$$
Also (\ref{eqn:detsig}) implies that the determinants of $K_k$ and $K_{k-1}$
are congruent modulo 4.  It follows that $\det Q_k\equiv \det K_k\pmod4$. Since
$\det Q_k$ and $\det K_k$ are equal up to sign and odd, $\det Q_k$ must be positive.  

On the other hand if $c_k$ is a
negative crossing then 
$$\det Q_k \equiv \det Q_{k-1}+2,\ \det K_k \equiv \det K_{k-1}+2\pmod 4,$$
and we again find $\det Q_k$ to be positive.  By induction $\det Q_k$ is positive for all $k$.\endpf

Finally note that we may reorient any of the link components
$C_1,\dots,C_{p+n}$ without changing the resulting Dehn surgery.
Also by rational handlesliding as in Proposition \ref{prop:slide}
we may change the linking matrix by ``adding'' $\pm2C_j$ to $C_i$
for any $i,j$.  These operations preserve the congruence classes
modulo 4 of the diagonal. The last claim in the statement of
Theorem \ref{thm:mainthm} now follows from the following lemma.

\begin{lemma}
\label{lem:mod2}
Any matrix $P\in GL(r,\zz)$ which is congruent to the identity modulo 2 may be obtained
from the identity by a sequence of row operations, each of which is either multiplying a row
by $-1$ or adding an even multiple of one row to another.
\end{lemma}
\proof
Let ${\bf b}=(b_1,\ldots,b_r)$ be an element of $\zz^r$
with $\gcd(b_1,\ldots,b_r)=1$.  Assume $b_i\ge0$ for all $i$, and that $b_1$ is
odd but the other components $b_2,\ldots,b_r$ are even.
Let $b_j$ be the least positive component.
By subtracting even multiples of $b_j$ and then possibly changing
sign, we may replace every other component $b_i$ by $b_i'$, with $0\le b_i'\le b_j$.  By the $\gcd$
condition, the least positive $b_i'$ is less than $b_j$ unless $b_j=j=1$.  By iterating this
procedure we see that $\bf{b}$ may be reduced to $(1,0,\ldots,0)$.

Now suppose $P\in GL(r,\zz)$ is congruent to $I$ modulo 2, and let
$\bf{b}$ be the first column of $Q$. The argument just given shows
that $Q$ may be replaced by a matrix with $(1,0,\ldots,0)$ in the
first column using the specified row operations.  Then replacing
the second column with $(*,1,0,\ldots,0)$ by row operations on the
last $r-1$ rows, and so on, we see that we may reduce $P$ to $I$
in this manner.
\endpf


\section{Heegaard Floer homology}
\label{sec:os}

In this section we recall some properties of the Heegaard Floer homology
invariants of \ozsvath\ and \szabo.  Details are to be found in their papers,
in particular \cite{os4,os6,osu1}.

Let $Y$ be an oriented rational homology three-sphere.  Recall that the space
$\spinc(Y)$ of \spinc structures on $Y$ is isomorphic to $H^2(Y;\zz)$.  If
$|H^2(Y;\zz)|$ is odd then there is a canonical isomorphism which takes the unique
spin structure to zero; this gives $\spinc(Y)$ a group structure.

Fixing a \spinc structure $\spincs$, the Heegaard Floer homology $HF^+(Y;\spincs)$
is a $\qq$-graded abelian group with an action by $\zz[U]$,
where $U$ lowers the grading by 2.  The \emph{correction term} invariant
is a rational number $d(Y,\spincs)$; it is defined to be the lowest grading of
a nonzero homogeneous element of $HF^+(Y;\spincs)$
which is in the image of $U^n$ for all  $n\in\nn.$
These have the property that $d(Y,\spincs)=-d(-Y,\spincs)$, where
$-Y$ denotes $Y$ with the opposite orientation.
We will describe below how these correction terms
may be computed in certain cases.

Now let $X$ be a positive-definite four-manifold with boundary $Y$.  Then it
is shown in \cite{os4} that for any \spinc structure $\spincs$ on $X$,
\begin{eqnarray}
c_1(\spincs)^2-b_2(X)&\ge& 4d(Y,\spincs|_Y),\label{eqn:ineq}\\
\mbox{and}\quad
c_1(\spincs)^2-b_2(X)&\equiv& 4d(Y,\spincs|_Y) \pmod2.\label{eqn:cong}
\end{eqnarray}
This means that the correction terms of $Y$ may be used to give an obstruction
to $Y$ bounding a four-manifold $X$ with a given positive-definite intersection form.
We will now elaborate on how this may be checked in practice.

Suppose for simplicity that $X$ is simply-connected and that $|H^2(Y;\zz)|$ is odd.
Let $r$ denote the second Betti number of $X$.
Fix a basis for $H_2(X;\zz)$ and thus an isomorphism
$$H_2(X;\zz)\cong\zz^r.$$
Let $Q$ be the matrix of the intersection pairing of $X$ in this basis; thus $Q$
is a symmetric positive-definite $r\times r$ integer matrix with $\det Q=|H^2(Y;\zz)|$.
The dual basis gives an isomorphism between
the second cohomology $H^2(X;\zz)$ and $\zz^r$.
The set $\{c_1(\spincs)\,\,|\,\,\spincs\in\spinc(X)\}\subset H^2(X;\zz)$
of first Chern classes
of \spinc structures is equal to the set of characteristic covectors
$Char(Q)$
for $Q$.  These in turn are elements $\xi$ of $\zz^r$ whose components
$\xi_i$ are congruent modulo 2 to the corresponding diagonal entries $Q_{ii}$ of $Q$.
The square of the first
Chern class of a \spinc structure is computed using the pairing induced by $Q$
on $H^2(X;\zz)$; in our choice of basis this is given by $\xi^T Q^{-1}\xi$.

The long exact sequence of the pair $(X,Y)$ yields the following short exact sequence:
$$0\longrightarrow\zz^r\stackrel{Q}{\longrightarrow}\zz^r
\longrightarrow H^2(Y;\zz)\longrightarrow0.$$
As in the introduction, define a function
$$m_Q:\zz^r/Q(\zz^r)\to\qq$$
by
$$m_Q(g) = \min\left\{\left.\frac{\xi^T Q^{-1}\xi-r}4\,\,\right|
\,\,\xi\in Char(Q),\,[\xi]=g\right\}.$$

In computing $m_Q$ it suffices to consider
characteristic covectors $\xi=(\xi_1,\dots,\xi_r)$ whose components are smaller
in absolute value than the corresponding diagonal entries of $Q$:
$$-Q_{ii}\le\xi_i\le Q_{ii}.$$
If, say, $\xi_i>Q_{ii}$, subtract twice the $i$th column of $Q$ from $\xi$ to
see that $\xi^TQ^{-1}\xi$ is not minimal.
A more difficult argument in \cite{os6} shows that it suffices to restrict to
$$-Q_{ii}\le\xi_i\le Q_{ii}-2.$$
Thus it is straightforward, if tedious, to compute $m_Q$ for a given
positive-definite matrix $Q$.

The conditions (\ref{eqn:ineq}) and (\ref{eqn:cong}) may now be expressed as follows:

\begin{theorem}[\ozsvath-\szabo]
\label{thm:os} Let $Y$ be a rational homology three-sphere which
is the boundary of a simply-connected positive-definite
four-manifold $X$, with $|H^2(Y;\zz)|$ odd.  If the intersection pairing of $X$ is
represented in a basis by the matrix $Q$ then there exists a group
isomorphism
$$\phi:\zz^r/Q(\zz^r)\to\spinc(Y)$$
with
\begin{eqnarray}
m_Q(g)&\ge&d(Y,\phi(g)),\label{eqn:ineqQ}\\
\mbox{and}\quad
m_Q(g)&\equiv&d(Y,\phi(g)) \pmod2\label{eqn:congQ}
\end{eqnarray}
for all $g\in\zz^r/Q(\zz^r)$.
\end{theorem}

The four-manifold $X$ is said to be \emph{sharp} if equality holds in (\ref{eqn:ineqQ}).
In this case the correction terms for $Y$
can be computed using the function $m_Q$ described above.
Also, if a rational homology sphere $Y$ bounds a negative-definite four-manifold $X$
such that $-X$ is sharp, then the correction terms for $Y$ can be computed
using the formula $d(Y,\spincs)=-d(-Y,\spincs)$.  Note
that if $K$ is a knot in $S^3$ then the standard orientation on $S^3$ induces an
orientation on $\Sigma(K)$; letting $r(K)$ denote the reflection of $K$, we have
$\Sigma(r(K))\cong-\Sigma(K)$.

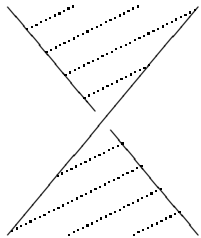
\begin{figure}[htbp] 
\begin{center}
\ifpic
\leavevmode
\begin{xy}
0;/r6pc/:
(0,0)*{}="1";
(1,0)*{}="2";
(1,1.2)*{}="3";
(0,1.2)*{}="4";
"1";"3" **\crv{}; \POS?(.5)*{\hole}="x"; 
"2";"x" **\crv{}; 
"x";"4" **\crv{}; 
(0.1,1.08);(0.34,1.2) **\dir{.};
(0.2,0.96);(0.68,1.2) **\dir{.};
(0.3,0.84);(0.98,1.18) **\dir{.};
(0.4,0.72);(0.74,0.89) **\dir{.};
(0.6,0.48);(0.26,0.31) **\dir{.};
(0.7,0.36);(0.02,0.02) **\dir{.};
(0.8,0.24);(0.32,0) **\dir{.};
(0.9,0.12);(0.66,0) **\dir{.};
\end{xy}
\else \vskip 5cm \fi
\begin{narrow}{0.3in}{0.3in}
\caption{
\bf{Colouring convention for alternating knot diagrams.}}
\label{fig:colour}
\end{narrow}
\end{center}
\end{figure}

In particular let $K$ be an alternating knot with double branched cover
$\Sigma(K)$.  Let $G$ denote the
positive-definite Goeritz matrix computed from an alternating diagram for $K$ as follows.
Colour the knot diagram in chessboard fashion according
to the convention shown in Figure \ref{fig:colour}.  (Note that this is the opposite
convention to that used in \cite{osu1}, since they use the negative-definite
Goeritz matrix.)  Let $v_1,\ldots,v_{k+1}$ denote the vertices of the white graph.
Then $G$ is the
$k\times k$ symmetric matrix $(g_{ij})$ with entries
$$g_{ij}=\left\{\begin{array}{ll}
\mbox{the number of edges containing $v_i$}&\mbox{if $i=j$}\\
\mbox{minus the number of edges joining $v_i$ and $v_j$ }&\mbox{if $i\ne j$}\\
\end{array}\right.$$
for $i,j=1,\ldots,k$.
It is shown in \cite[Proposition 3.2]{osu1} that $G$
represents the intersection pairing of a sharp four-manifold bounded by $\Sigma(K)$.
Thus the correction terms for $\Sigma(K)$ are given by $m_G$ (for any choice of
alternating diagram and any ordering of the white regions).  Also it follows from
\cite{gl} that with this colouring convention, the signature of $K$ is given by
$$\sigma(K) = k - \mu,$$
where $\mu$ is the number of positive crossings in the alternating diagram used to
compute $G$.

Also if $K$ is a Montesinos knot then the double branched cover $\Sigma(K)$ is
a Seifert fibred space which is given as the boundary of a plumbing of
disk bundles over $S^2$.  This plumbing is determined (nonuniquely)
by the Montesinos invariants
which specify $K$.  After possibly reflecting $K$ we may choose the plumbing
so that its intersection pairing
is represented by a positive-definite matrix $P$.
It is shown in \cite{os6} that the plumbing is sharp,
so that the correction terms for $\Sigma(K)$ are given by $m_P$.
(See \cite{bs} for a description of Montesinos knots and their
branched double covers.)

\begin{remark}
Checking the congruence condition (\ref{eqn:cong}) alone is equivalent to
checking that the intersection pairing of $X$ presents the linking pairing
of $Y$; see \cite{mac} for a detailed discussion.
\end{remark}


\section{Obstruction to unknotting}
\label{sec:u}

In this section we prove Theorems \ref{thm:u2} and \ref{thm:u}.

Let $\calq(r,\delta)$ denote the set of positive-definite symmetric integer matrices of
rank $r$ and determinant $\delta$, on which $GL(r,\zz)$ acts by $P\cdot Q=PQP^T$
with finite quotient (see e.g. \cite{cassels}).
Let $\calq(r,\delta)_2\subset\calq(r,\delta)$ (resp. $GL(r,\zz)_2\subset GL(r,\zz)$)
denote the subset (resp. subgroup) consisting
of matrices which are congruent to the identity modulo 2.
Then the subset $\calq(r,\delta)_2/GL(r,\zz)$ is clearly finite, and thus so is
$\calq(r,\delta)_2/GL(r,\zz)_2$ since $GL(r,\zz)_2$ is a finite index subgroup of $GL(r,\zz)$.

 \

\noindent{\bf Proof of Theorem \ref{thm:u}}.
By Theorem \ref{thm:mainthm}, the unknotting hypothesis implies that $\Sigma(K)$ is given by
Dehn surgery on a link in $S^3$ with linking matrix $\frac12Q_i$ for some $i$,
where $n$ of the diagonal entries of $Q_i$ are congruent to 3 modulo 4.  By Lemma
\ref{lem:links}, $\Sigma(K)$ bounds the 2-handlebody $W$ specified by an integer-framed link
with positive-definite linking matrix $\widetilde{Q}_i$, which then represents the
intersection pairing of
$W$.  The conclusion now follows from Theorem \ref{thm:os}.\endpf

 \

\noindent{\bf Proof of Theorem \ref{thm:u2}}.
Theorem \ref{thm:u2} follows from Theorem \ref{thm:u} since a finite set of representatives of
$\calq(2,\delta)_2/GL(2,\zz)_2$ is given by the set of matrices
$$\left\{Q=\left(\left.\begin{matrix}2m_1-1 & 2a\\ 2a & 2m_2-1 \end{matrix}\right)\,\,\right|\,\,
\det Q=\delta,\, 0\le a<m_1\le m_2\right\},$$
and since the correction terms $d(\Sigma(K),\spincs)$ may be computed using a
positive-definite Goeritz matrix $G$ when $K$ is alternating.\endpf

\begin{remark}
Theorems \ref{thm:u2} and \ref{thm:u} do not use all
of the information from Theorem \ref{thm:mainthm}.  We have only used the information
about the intersection pairing of the four-manifold $W$ bounded by $\Sigma(K)$, and not
the fact that $W$ is a surgery cobordism arising from a half-integral surgery.
Comparing to Theorem 1.1 in \cite{osu1}, we have generalised conditions (1) and (2) to the
case of $u(K)>1$ but not the symmetry condition (3).  It is to be hoped that the symmetry
condition may also be generalised in some way, leading to a stronger obstruction
and computation of some more unknotting numbers.
\end{remark}


\section{Examples}
\label{sec:knots}

\noindent{\bf Proof of Corollary \ref{cor:u}.}
For each knot in Corollary \ref{cor:u}
we distinguish between $K$ and its reflection $r(K)$ by specifying that
$K$ has positive signature.

\begin{figure}[htbp]
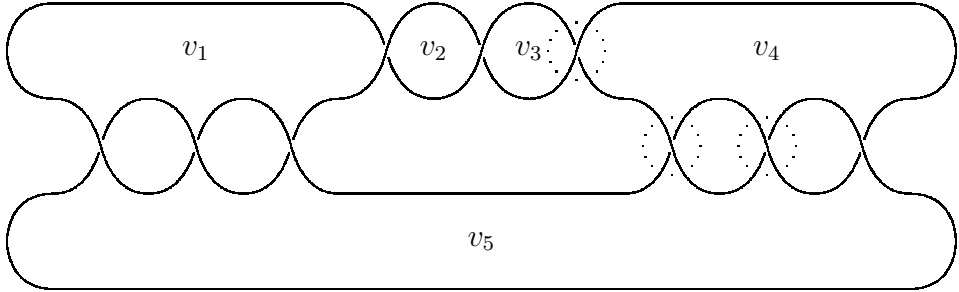

\begin{center}
\ifpic
\leavevmode
\xygraph{
!{0;/r3.0pc/:}
!{\hcap[-1]}[dd]
!{\hcap[-1]}[uu]
!{\xcaph[1]@(0)}[ld]
!{\htwist}[ldd]
!{\xcaph[1]@(0)}[uuu]
!{\xcaph[1]@(0)}[ld]
!{\htwist}[ldd]
!{\xcaph[1]@(0)}[uuu]
!{\xcaph[1]@(0)}[ld]
!{\htwist}[ldd]
!{\xcaph[1]@(0)}[uuu]
!{\htwistneg}[ldd]
!{\xcaph[1]@(0)}[ld]
!{\xcaph[1]@(0)}[uuu]
!{\htwistneg}[ldd]
!{\xcaph[1]@(0)}[ld]
!{\xcaph[1]@(0)}[uuu]
!{\htwistneg}[ldd]
!{\xcaph[1]@(0)}[ld]
!{\xcaph[1]@(0)}[uuu]
!{\xcaph[1]@(0)}[ld]
!{\htwist}[ldd]
!{\xcaph[1]@(0)}[uuu]
!{\xcaph[1]@(0)}[ld]
!{\htwist}[ldd]
!{\xcaph[1]@(0)}[uuu]
!{\xcaph[1]@(0)}[ld]
!{\htwist}[ldd]
!{\xcaph[1]@(0)}[uuu]
!{\hcap}[dd]
!{\hcap}[u(0.5)l(1.5)]
!{\ellipse(0.3){.}}[l]
!{\ellipse(0.3){.}}[ul]
!{\ellipse(0.3){.}}[l(4)]
!{*{v_1}}[r(2.5)]
!{*{v_2}}[r]
!{*{v_3}}[r(2.5)]
!{*{v_4}}[l(3)d(2)]
!{*{v_5}}
}
\else \vskip 5cm \fi
\begin{narrow}{0.3in}{0.3in}
\caption{
{\bf{The knot $9_{10}=S(33,23)$.}}
Note that changing the circled crossings will give the unknot.  The labels
$v_1,\ldots, v_5$ correspond to vertices of the white graph.}
\label{fig:9_10}
\end{narrow}
\end{center}
\end{figure}

We start with the knot $K=9_{10}$ shown in Figure \ref{fig:9_10}.  This is the two-bridge
knot $S(33,23)$.  It has signature $4$, and it is easy to see that 3 crossing changes
suffice to unknot it.  Thus the unknotting number is either 2 or 3, and if it can
be unknotted by changing two crossings then both are negative ($p=0$ and $n=2$).

With the white regions labelled as shown in the figure, the
Goeritz matrix is
$$G=\left(\begin{matrix}
4 & -1 & 0 & 0 \\
-1 & 2 & -1 & 0 \\
0 & -1 & 2 & -1 \\
0 & 0 & -1 & 4 \\
\end{matrix}\right).$$

Using $m_G$, we find the correction terms of $\Sigma(K)$ to be:

$$A=\left\{
\begin{array}{rrrrrrrrrrr}
-1,& -\frac{23}{33},& \frac{7}{33},& -\frac{3}{11},& -\frac{5}{33},& \frac{19}{33},
& -\frac{1}{11},& -\frac{5}{33},& \frac{13}{33},& -\frac{5}{11},& -\frac{23}{33},\\
&&&&&&&&&&\\
-\frac{1}{3},& \frac{7}{11},& \frac{7}{33},& \frac{13}{33},& \frac{13}{11},
& \frac{19}{33},& \frac{19}{33},& \frac{13}{11},& \frac{13}{33},& \frac{7}{33},
& \frac{7}{11},\\
&&&&&&&&&&\\
-\frac{1}{3},& -\frac{23}{33},& -\frac{5}{11},& \frac{13}{33},& -\frac{5}{33},
& -\frac{1}{11},& \frac{19}{33},& -\frac{5}{33},& -\frac{3}{11},& \frac{7}{33},
& -\frac{23}{33}\\
\end{array}
\right\}.$$
The order of this list corresponds to the cyclic group structure of
$\spinc(\Sigma(K))\cong H^2(\Sigma(K);\zz)$,
and the first element is the correction term of the spin structure.

The determinant of $9_{10}$ is $33$.  To find a matrix ${\widetilde{Q}}$
as in Theorem \ref{thm:u2} we need to find
$(m_1,a,m_2)$ with
$$(2m_1-1)(2m_2-1)-4a^2=33,$$
$$ 0\le a<m_1\le m_2,$$
and $m_1$ and $m_2$ are even.  There are two solutions: $(2,0,6)$ and $(4,2,4)$.
Computing $m_{\widetilde{Q}}$ for each of the matrices
$${\widetilde{Q}}_1=\left(\begin{matrix}
2 & 1 & 0 & 0\\
1 & 2 & 0 & 0\\
0 & 0 & 6 & 1\\
0 & 0 & 1 & 2
\end{matrix}\right),\quad
{\widetilde{Q}}_2=\left(\begin{matrix}
4 & 1 & 2 & 0\\
1 & 2 & 0 & 0\\
2 & 0 & 4 & 1\\
0 & 0 & 1 & 2
\end{matrix}\right)$$
yields the following lists:
$$B_1=\left\{
\begin{array}{rrrrrrrrrrr}
-1,& -\frac{5}{33},& \frac{13}{33},& \frac{7}{11},& \frac{19}{33},& \frac{7}{33},
& -\frac{5}{11},& \frac{19}{33},& \frac{43}{33},& -\frac{3}{11},& -\frac{5}{33},\\
&&&&&&&&&&\\
-\frac{1}{3},& -\frac{9}{11},& \frac{13}{33},& \frac{43}{33},& -\frac{1}{11},
& \frac{7}{33},& \frac{7}{33},& -\frac{1}{11},& \frac{43}{33},& \frac{13}{33},
& -\frac{9}{11},\\
&&&&&&&&&&\\
-\frac{1}{3},& -\frac{5}{33},& -\frac{3}{11},& \frac{43}{33},& \frac{19}{33},
& -\frac{5}{11},& \frac{7}{33},& \frac{19}{33},& \frac{7}{11},& \frac{13}{33},
& -\frac{5}{33}\\
\end{array}
\right\},$$
$$B_2=\left\{
\begin{array}{rrrrrrrrrrr}
-1,& -\frac{19}{33},& \frac{23}{33},& \frac{9}{11},& -\frac{7}{33},& -\frac{13}{33},
& \frac{3}{11},& -\frac{7}{33},& \frac{5}{33},& -\frac{7}{11},& -\frac{19}{33},\\
&&&&&&&&&&\\
\frac{1}{3},& \frac{1}{11},& \frac{23}{33},& \frac{5}{33},& \frac{5}{11},
& -\frac{13}{33},& -\frac{13}{33},& \frac{5}{11},& \frac{5}{33},& \frac{23}{33},
& \frac{1}{11},\\
&&&&&&&&&&\\
\frac{1}{3},& -\frac{19}{33},& -\frac{7}{11},& \frac{5}{33},& -\frac{7}{33},
& \frac{3}{11},& -\frac{13}{33},& -\frac{7}{33},& \frac{9}{11},& \frac{23}{33},
& -\frac{19}{33}\\
\end{array}
\right\}.$$

We claim that for both ${\widetilde{Q}}_1$ and ${\widetilde{Q}}_2$
it is impossible to find a group automorphism
$\phi$ of $\zz/33$ satisfying the required inequality and congruence conditions.
This is immediate in either case by considering the minimal elements
(excluding $-1$ which appears in all 3 lists).  We have the
entry $-9/11$ in $B_1$.  By inspection there is no element in $A$ which is less than
or equal to $-9/11$, and differs from it by a multiple of 2.
The same applies to $-7/11$ in $B_2$.  We conclude that $9_{10}$ cannot be unknotted
by two crossing changes and $u(9_{10})=3$.

Similar calculations show that $9_{13}, 9_{38}, 10_{53}, 10_{101}$ and $10_{120}$
cannot be unknotted with two crossing changes.
All of these knots are alternating, have signature four and cyclic $H^2(\Sigma(K);\zz)$.
By inspection of their diagrams (see e.g. \cite{knotinfo}),
all can be unknotted with three crossing changes.  For some details of the
calculations for these knots, see Table \ref{table:u}.  Note that we use
the knot diagrams from \cite{knotinfo} to compute the Goeritz matrices
for these knots, after possibly
reflecting to ensure positive signature.

\begin{figure}[htbp]
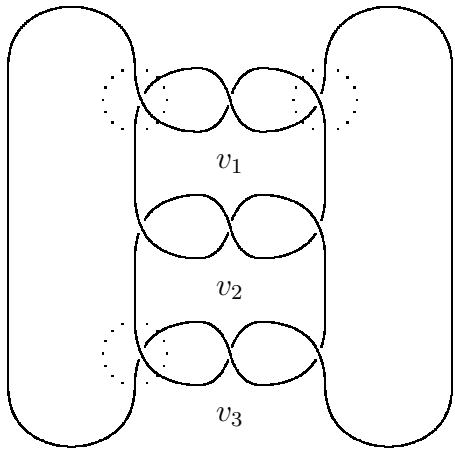

\begin{center}
\ifpic \leavevmode \xygraph{ !{0;/r2pc/:}
[r(3)]
!{\vcap[2]}
!{\xcapv[5]@(0)}[dddd]
!{\vcap[-2]}[u(5)r(5)]
!{\vcap[2]}[rr]
!{\xcapv[5]@(0)}[d(4)l(2)]
!{\vcap[-2]}[u(5)l(3)]
!{\hunder}
!{\htwist}
!{\hunder-}[dlll]
!{\xcapv@(0)}[r(3)u(1)]
!{\xcapv@(0)}[l(3)]
!{\hunder}
!{\htwist}
!{\hunder-}[dlll]
!{\xcapv@(0)}[rrru]
!{\xcapv@(0)}[lll]
!{\hunder}
!{\htwist}
!{\hunder-}[d(0.5)lll]
!{\ellipse(0.5){.}}[u(4)]
!{\ellipse(0.5){.}}[r(3)]
!{\ellipse(0.5){.}}[l(1.5)d]
!{*{v_1}}[d(2)]
!{*{v_2}}[d(2)]
!{*{v_3}}
}
\else \vskip 5cm \fi
\begin{narrow}{0.3in}{0.3in}
\caption{
{\bf{The Montesinos knot $9_{35} = M(0;(3,1),(3,1),(3,1))$.}}
}
\label{fig:9_35}
\end{narrow}
\end{center}
\end{figure}

Finally consider $K=9_{35}$, pictured in Figure \ref{fig:9_35}.
It has signature 2 and can be unknotted with 3 crossing changes.
The Goeritz matrix from the figure is
$$G=\left(\begin{matrix}
6 & -3\\
-3 & 6
\end{matrix}\right).$$
We note that this presents $H^2(\Sigma(K);\zz)$ which is thus
2-cyclic; this shows (by Montesinos' theorem for example but by an
inequality originally due to Wendt) that $u(K)\ge 2$.  We can use
$m_G$ to compute the correction terms of $\Sigma(K)$, which are
$$A=\left[
\begin{matrix}
-\frac{1}{2} & \frac{19}{18} & -\frac{5}{18} & \frac{3}{2} & \frac{7}{18} &
\frac{7}{18} & \frac{3}{2} & -\frac{5}{18} & \frac{19}{18}\\
&&&&&&&&\\
\frac{1}{6} & -\frac{5}{18} & \frac{7}{18} & \frac{1}{6} & \frac{19}{18} &
\frac{19}{18} & \frac{1}{6} & \frac{7}{18} & -\frac{5}{18}\\
&&&&&&&&\\
\frac{1}{6} & -\frac{5}{18} & \frac{7}{18} & \frac{1}{6} & \frac{19}{18} &
\frac{19}{18} & \frac{1}{6} & \frac{7}{18} & -\frac{5}{18}
\end{matrix}
\right].$$
Here the rectangular array shows the $\zz/3\oplus\zz/9$ group structure; the top left
entry is the correction term of the spin structure.

Suppose that $9_{35}$ may be unknotted by changing one positive and one negative
crossing.  The only matrix which satisfies the conditions of Theorem \ref{thm:u2}
and which presents $\zz/3\oplus\zz/9$ is
$${\widetilde{Q}}=\left(
\begin{matrix}
2 & 1 & 0 & 0\\
1 & 2 & 0 & 0\\
0 & 0 & 5 & 1\\
0 & 0 & 1 & 2\\
\end{matrix}
\right).$$
Computing $m_{\widetilde{Q}}$
yields another array whose minimal entry is $-17/18$; we conclude that there
is no automorphism $\phi$ of $\zz/3\oplus\zz/9$ satisfying
the conclusion of Theorem \ref{thm:u2}.

This is not enough to rule out the possibility that $u(9_{35})=2$; it does however
show that if $9_{35}$ can be unknotted by two crossing changes, then they are both
negative crossings.  Using the value of the Jones polynomial at $e^{i\pi/3}$,
Traczyk has shown in \cite{t} that if $9_{35}$ can be unknotted by changing two
crossings, then the crossings have different signs.  We conclude that $u(9_{35})=3$.
\endpf

\begin{figure}[htbp]
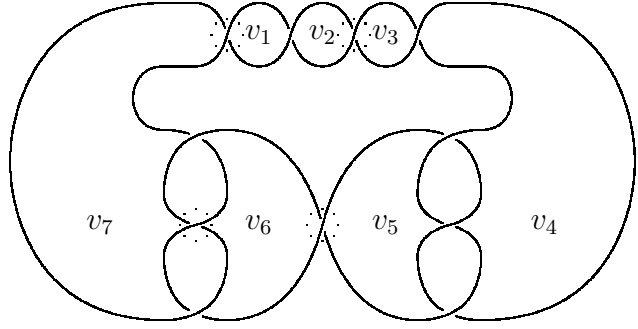

\begin{center}
\ifpic
\leavevmode
\xygraph{
!{0;/r2pc/:}
!{\hcap[-1]}
[u]!{\hcap[-5]}
[d(4)]!{\vunder[-1]}
[u(2)]!{\vtwist}
[u(2)]!{\vunder}
[ur]!{\htwistneg[3]}
[u(2)l(2)]
!{\xcaph[0.5]@(0)}
[dl]!{\xcaph[0.5]@(0)}
[u(1)l(0.5)]!{\htwistneg}
!{\htwistneg}
!{\htwistneg}
!{\htwistneg}
!{\xcaph[0.5]@(0)}
[dl]!{\xcaph[0.5]@(0)}
[l(0.5)]!{\hcap}
[u]!{\hcap[5]}
[d(4)l(1)]!{\vunder[-1]}
[u(2)]!{\vtwist}
[u(2)]!{\vunder}
[u(2.5)l(1)]!{\ellipse(0.25){.}}
[l(2)]!{\ellipse(0.25){.}}
[l(0.5)d(3)]!{\ellipse(0.25){.}}
[r(2)]!{\ellipse(0.25){.}}
[l(3.5)]!{*{v_7}}
[u(3)r(2.5)]!{*{v_1}}
[r]!{*{v_2}}
[r]!{*{v_3}}
[d(3)r(2.5)]!{*{v_4}}
[l(2.5)]!{*{v_5}}
[l(2)]!{*{v_6}}
}
\else \vskip 5cm \fi
\begin{narrow}{0.3in}{0.3in}
\caption{
{\bf{The two-bridge knot $S(51,35)$, or $11a365$.}}
}
\label{fig:11a365}
\end{narrow}
\end{center}
\end{figure}

\noindent{\bf Proof of Corollary \ref{cor:11a365}.}
The two-bridge knot $K=S(51,35)$ is listed in \cite{knotinfo} as $11a365$ and is shown in
Figure \ref{fig:11a365}.  It has signature $6$, and from the diagram we see that it may be
unknotted by changing 4 crossings.  We will apply Theorem \ref{thm:u} to show that it does
not have $u(K)=n=3$.

Note that $\det K=51$.
In order to apply Theorem \ref{thm:u} we first need to find a set of representatives of the
finite quotient $\calq(3,51)_2/GL(3,\zz)_2$ 
with all diagonal entries conjugate to 3 modulo 4.  According
to \cite{j}, a complete set of representatives of $\calq(3,51)/GL(3,\zz)$ is
given by the (Eisenstein reduced) matrices
$$\left(\begin{matrix}1&0&0\\0&1&0\\0&0&51\end{matrix}\right),
\left(\begin{matrix}1&0&0\\0&2&1\\0&1&26\end{matrix}\right),
\left(\begin{matrix}1&0&0\\0&3&0\\0&0&17\end{matrix}\right),
\left(\begin{matrix}1&0&0\\0&4&1\\0&1&13\end{matrix}\right),
\left(\begin{matrix}1&0&0\\0&5&2\\0&2&11\end{matrix}\right),
\left(\begin{matrix}1&0&0\\0&6&3\\0&3&10\end{matrix}\right),$$
$$\left(\begin{matrix}2&0&1\\0&3&0\\1&0&9\end{matrix}\right),
\left(\begin{matrix}2&1&0\\1&2&0\\0&0&17\end{matrix}\right),
\left(\begin{matrix}3&0&1\\0&3&0\\1&0&6\end{matrix}\right),
\left(\begin{matrix}3&1&1\\1&4&0\\1&0&5\end{matrix}\right),
\left(\begin{matrix}4&1&2\\1&4&2\\2&2&5\end{matrix}\right).$$
Note that if $P\in GL(3,\zz)$ satisfies $PP^T\equiv I\pmod2$, then $P$ is conjugate to a
permutation matrix modulo 2.  Thus if $P\in GL(3,\zz)$ and
$Q,PQP^T\in\calq(3,51)_2$ then $Q$ and
$PQP^T$ have the same number of diagonal entries conjugate to 3 modulo 4.
We therefore eliminate the forms represented by
$\left(\begin{matrix}1&0&0\\0&1&0\\0&0&51\end{matrix}\right)$,
$\left(\begin{matrix}1&0&0\\0&3&0\\0&0&17\end{matrix}\right)$,
$\left(\begin{matrix}1&0&0\\0&5&2\\0&2&11\end{matrix}\right)$.
For each remaining form in the list, we look for a basis in which the form is congruent
to the identity modulo 2.  If no such basis exists, or if we find that some diagonal entry is
not conjugate to 3 modulo 4, we eliminate the form.  This leaves us with the following four
forms to consider:
$\left(\begin{matrix}1&0&0\\0&2&1\\0&1&26\end{matrix}\right)\sim
\left(\begin{matrix}3&2&0\\2&27&26\\0&26&27\end{matrix}\right)$,
$\left(\begin{matrix}1&0&0\\0&6&3\\0&3&10\end{matrix}\right)\sim
\left(\begin{matrix}11&4&-6\\4&7&4\\-6&4&11\end{matrix}\right)$,
$\left(\begin{matrix}2&1&0\\1&2&0\\0&0&17\end{matrix}\right)\sim
\left(\begin{matrix}19&18&18\\18&19&16\\18&16&19\end{matrix}\right)$, and
$\left(\begin{matrix}3&0&1\\0&3&0\\1&0&6\end{matrix}\right)\sim
\left(\begin{matrix}3&0&-2\\0&3&0\\-2&0&7\end{matrix}\right)$.

From Figure \ref{fig:11a365} we may write down the Goeritz matrix $G$; the
correction terms $\{d(\Sigma(K),\spincs)\}$ are then given by $m_G$.
For each $Q$ in
$$\left\{
\left(\begin{matrix}3&2&0\\2&27&26\\0&26&27\end{matrix}\right),
\left(\begin{matrix}11&4&-6\\4&7&4\\-6&4&11\end{matrix}\right),
\left(\begin{matrix}19&18&18\\18&19&16\\18&16&19\end{matrix}\right),
\left(\begin{matrix}3&0&-2\\0&3&0\\-2&0&7\end{matrix}\right)\right\},$$
one may check that there is no isomorphism
$$\phi:\Gamma_{\widetilde{Q}}\to\spinc(\Sigma(K))$$
satisfying the conclusion of Theorem \ref{thm:u}.
We conclude that the unknotting number of $K$ is 4.\endpf

\begin{remark}
\label{rmk:reduce}
In the last step of the proof of Corollary \ref{cor:11a365} it is much quicker in some
cases to change basis before computing $m_{\widetilde{Q}}$, so as to work with a matrix
with smaller diagonal entries.  For example with
$Q=\left(\begin{matrix}19&18&18\\18&19&16\\18&16&19\end{matrix}\right)$,
we have
$$\widetilde{Q}=
\left(\begin{matrix}
10&1&9&0&9&0\\
1&2&0&0&0&0\\
9&0&10&1&8&0\\
0&0&1&2&0&0\\
9&0&8&0&10&1\\
0&0&0&0&1&2\end{matrix}\right)\sim
\left(\begin{matrix}
10&1&-1&0&-1&0\\
1&2&-1&0&-1&0\\
-1&-1&2&1&0&0\\
0&0&1&2&0&0\\
-1&-1&0&0&2&1\\
0&0&0&0&1&2\end{matrix}\right),$$
by subtracting the first basis vector from the third and fifth.
As a result we need to consider $2^5\cdot10$ characteristic covectors to compute
$m_{\widetilde{Q}}$ instead of $2^3\cdot10^3$.
\end{remark}

\begin{table}[!ht]
\begin{center}
\begin{tabular}{|c|c|c|c|c|c|} \hline
Knot & Goeritz matrix & $\ds\min_{g\ne0}\{m_G(g)\}$ & $(m_1,a,m_2)$ &
$\ds\min_{g\ne0}\{m_Q(g)\}$ \\
\hline
$9_{13}$ &
$\left(
\begin{matrix} 2&-1&0&0\\ -1&2&-1&0\\ 0&-1&4&-1\\ 0&0&-1&4 \end{matrix}
\right)$
& $-\frac{27}{37}$ & $(10,9,10)$ & $-\frac{33}{37}$\\
&&&&\\
$9_{38}$ &
$\left(
\begin{matrix} 4&-1&-1&0\\ -1&4&-2&0\\ -1&-2&4&-1\\ 0&0&-1&2 \end{matrix}
\right)$
& $-\frac{37}{57}$ & $(2,0,10)$ & $-\frac{51}{57}$\\
&&& $(6,4,6)$ & $-\frac{45}{57}$\\
&&&&\\
$10_{53}$ &
$\left(
\begin{matrix} 4&-1&0&0\\ -1&4&-1&-1\\ 0&-1&4&-1\\ 0&-1&-1&2 \end{matrix}
\right)$
& $-\frac{53}{73}$ & $(4,1,6)$ & $-\frac{59}{73}$\\
&&&&\\
$10_{101}$ &
$\left(
\begin{matrix} 2&-1&0&0\\ -1&4&-1&-1\\ 0&-1&4&-1\\ 0&-1&-1&4 \end{matrix}
\right)$
& $-\frac{59}{85}$ & $(6,3,6)$ & $-\frac{65}{85}$\\
&&& $(22,21,22)$ & $-\frac{81}{85}$\\
&&&&\\
$10_{120}$ &
$\left(
\begin{matrix} 4&-2&0&-1\\ -2&4&-1&0\\ 0&-1&4&-2\\ -1&0&-2&4 \end{matrix}
\right)$
& $-\frac{69}{105}$ & $(2,0,18)$ & $-\frac{99}{105}$\\
&&& $(4,0,8)$ & $-\frac{91}{105}$\\
&&&&\\
&&& $(6,2,6)$ & $-\frac{83}{105}$\\
&&&&\\
&&& $(10,8,10)$ & $-\frac{93}{105}$\\
&&&&\\
\hline
\end{tabular}
\vskip5mm
\begin{narrow}{0.3in}{0.3in}
\caption{
{\bf{Data for knots in Corollary \ref{cor:u}.}}
The fourth column contains possible coefficients of the matrix
$Q$ in Theorem \ref{thm:u2}.
}
\label{table:u}
\end{narrow}
\end{center}
\end{table}


\end{document}